\documentclass[a4paper,11pt]{article}
\usepackage{graphicx}
\usepackage{amsmath,amsthm,amssymb,enumerate}
\usepackage{euscript,mathrsfs}
\usepackage{dsfont}
\usepackage{xcolor}
\usepackage{empheq}
\usepackage{geometry}

\usepackage{bm}
\usepackage[linkcolor=blue,colorlinks=true]{hyperref}
\catcode`\@=11 \@addtoreset{equation}{section}

\catcode`\@=12

\usepackage{stmaryrd}
\allowdisplaybreaks

\usepackage[labelformat=simple]{subcaption}

\usepackage{enumitem}
\usepackage{url}

\usepackage{tikz}
\usetikzlibrary{patterns}
\usepackage[normalem]{ulem}
\usepackage{cancel}

\newtheorem{Theorem}{Theorem}[section]
\newtheorem{Proposition}[Theorem]{Proposition}
\newtheorem{Lemma}[Theorem]{Lemma}
\newtheorem{Corollary}[Theorem]{Corollary}

\newtheorem{Definition}[Theorem]{Definition}

\newtheorem{Remark}[Theorem]{Remark}
\newcommand{\bRemark}[1]{
			\begin{Remark} \label{R#1} }
\newcommand{\eR}{\end{Remark}}

\newcommand{\ith}{i^{\rm th}}
\newcommand{\ha}{h^{\alpha}}

\newcommand{\RE}[2]{R_E\left((#1)\mid(#2)\right)}
\newcommand{\EH}[2]{\mathbb{E}_{\cal H}\big((#1)|(#2)  \big)}

\newcommand{\vrh}{\vr_h}
\newcommand{\vuh}{\bm{u}_h}
\newcommand{\vu}{\bm{u}}
\newcommand{\ve}{\bm{e}}

\newcommand{\auh}{\avs{\vuh}}
\newcommand{\vth}{\vt_h}
\newcommand{\vthout}{\vth^{\rm out}}

\newcommand{\uih}{u_{i,h}}
\newcommand{\ujh}{u_{j,h}}

\newcommand{\Dhuh}{\bD_h(\vuh)}

\newcommand{\Up}{{\rm Up}}

\newcommand{\Fup}{F_h^\eps}

\newcommand{\tor}{\mathbb{T}^d}

\newcommand{\bfv}{\bm{v}}

\newcommand{\bfphi}{\boldsymbol{\phi}}

\newcommand{\Piq}{\Pi_Q}
\newcommand{\Piw}{\Pi_W}
\newcommand{\Piwi}{\Pi_W^{(i)}}

\newcommand{\sumj}{\sum_{j=1}^d}

\newcommand{\ds}{\,{\rm d}S_x}

\newcommand{\grid}{{\cal T}_h}
\newcommand{\dgrid}{\mathcal{D}}
\newcommand{\Di}{\dgrid_i}

\newcommand{\TS}{\Delta t}

\newcommand{\Divh}{{\rm div}_h}
\newcommand{\Gradh}{\nabla_h}

\newcommand{\Gradedge}{\nabla_\faces}

\newcommand{\pdedgej}{\eth_ \faces^{(j)}}

\newcommand{\pdmeshj}{\eth_{\cal T}^{(j)}}

\newcommand{\co}[2]{{\rm co}\{ #1 , #2 \}}
\newcommand{\Ov}[1]{\overline{ #1 } }
\newcommand{\avs}[1]{\left\{\hspace{-3pt}\left\{ #1 \right\} \hspace{-3pt}\right\} }

\newcommand{\aleq}{\stackrel{<}{\sim}}

\newcommand{\Un}[1]{\underline{#1}}
\newcommand{\vr}{\varrho}

\newcommand{\vt}{\vartheta}
\newcommand{\vm}{\bm{m}}

\newcommand{\vn}{\bm{n}}

\newcommand{\tvr}{\widetilde \vr}
\newcommand{\tvu}{{\widetilde \vu}}
\newcommand{\tvt}{\widetilde \vt}
\newcommand{\tp}{{\widetilde p}}
\newcommand{\ts}{\widetilde s}

\newcommand{\Ovt}{\Ov{\vt}}
\newcommand{\Uvt}{\underline{\vt}}

\newcommand{\vT}{\tvt}

\newcommand{\vc}[1]{{\bm #1}}

\newcommand{\Div}{{\rm div}_x}
\newcommand{\Grad}{\nabla_x}

\newcommand{\dx}{\,{\rm d} {x}}

\newcommand{\dt}{\,{\rm d} t }
\newcommand{\dxdt}{\,{\rm d} {x}{\rm d} t }

\newcommand{\jump}[1]{\left\llbracket#1\right\rrbracket}
\newcommand{\abs}[1]{| #1|}
\newcommand{\Abs}[1]{\left| #1 \right|}

\newcommand{\norm}[1]{\left\lVert#1\right\rVert}

\newcommand{\intO}[1]{\int_{\tor} #1 \dx}
\newcommand{\intOB}[1]{\int_{\tor} \left(#1\right) \dx}
\newcommand{\intTd}[1]{\int_{\tor} #1 \dx}
\newcommand{\intTdB}[1]{\int_{\tor} \left(#1\right) \dx}
\newcommand{\intfacesint}[1]{\int_{\facesint}{ #1 \ds} }

\newcommand{\Of}{\tor}

\newcommand{\intSh}[1] {\int_{\sigma} #1 \ds }

\newcommand{\intK}[1] {\int_{K} #1 \dx }

\newcommand{\intn}{\int_{0}^{t^{n+1}}}
\newcommand{\intTO}[1]{ \int_0^\tau \int_{\tor} #1 \dxdt}
\newcommand{\intTOB}[1]{  \int_0^\tau \int_{\tor} \left( #1 \right) \dxdt}

\newcommand{\vv}{\bm{v}}

\newcommand{\R}{\mathbb{R}}
\newcommand{\I}{\mathbb{I}}

\newcommand{\intTor}[1]{\int_{\tor} #1 \ \dx}

\def\softd{{\leavevmode\setbox1=\hbox{d}%
          \hbox to 1.05\wd1{d\kern-0.4ex{\char039}\hss}}}

\newcommand{\bdot}{\boldsymbol{\cdot}}

\definecolor{Cgrey}{rgb}{0.85,0.85,0.85}
\definecolor{Cblue}{rgb}{0.50,0.85,0.85}
\definecolor{Cred}{rgb}{1,0,0}
\definecolor{fancy}{rgb}{0.10,0.85,0.10}
\definecolor{forestgreen}{rgb}{0.13, 0.55, 0.13}
\definecolor{pinegreen}{rgb}{0.0, 0.47, 0.44}

\newcommand{\cblue}{\color{blue}}

\date{}
\newcommand{\pd}{\partial}
\newcommand{\Hc}{\mathcal{H}_{\vT}}
\newcommand{\eps}{\varepsilon}
\newcommand{\faces}{\mathcal{E}}
\newcommand{\facesi}{\faces _i}
\newcommand{\facesK}{\faces(K)}

\newcommand{\facesint}{\faces}

\newcommand{\bQh}{ Q_h}

\newcommand{\bD}{\mathbb D}
\newcommand{\bS}{\mathbb S}
\newcommand{\tbS}{\widetilde{\bS}}
\newcommand{\bDu}{\bD(\vu)}

\newcommand{\muh}{h^\eps}

\newcommand{\Laph}{\Delta_h}
\newcommand{\Laphj}{\Delta_h^{(j)}}

\newcommand{\br}{ \nonumber \\ }

\allowdisplaybreaks

\begin{document}

\title{Error estimates of a finite volume method for the compressible Navier--Stokes--Fourier system}
\author{
Danica Basari\'c\thanks{D.B., H.M. and B.S. have received funding from
the Czech Sciences Foundation (GA\v CR), Grant Agreement 21--04211S.
The Institute of Mathematics of the Czech Academy of Sciences is supported by RVO:67985840.
}
\and M\'aria Luk\'a\v{c}ov\'a-Medvi\softd ov\'a\thanks{The work of M.L. and Y. Y. was supported by the Deutsche Forschungsgemeinschaft (DFG, German Research
Foundation) - Project number 233630050 - TRR 146 as well as by TRR 165 Waves to Weather.
M.L. is grateful
to the Gutenberg Research College and Mainz Institute of Multiscale Modelling for supporting her research. }
\and
Hana Mizerov\'a$^{*, \ddagger}$
\and Bangwei She$^{*,\spadesuit}$
\and Yuhuan Yuan$^\dagger$
}

\date{}

\maketitle

\centerline{$^*$ Institute of Mathematics of the Czech Academy of Sciences}
\centerline{\v Zitn\' a 25, CZ-115 67 Praha 1, Czech Republic}
\centerline{(basaric,mizerova,she)@math.cas.cz}

\bigskip
\centerline{$^\dagger$ Institute of Mathematics, Johannes Gutenberg-University Mainz}
\centerline{Staudingerweg 9, 55 128 Mainz, Germany}
\centerline{(lukacova,yuhuyuan)@uni-mainz.de}

\bigskip
\centerline{$^\ddagger$ Department of Mathematical Analysis and Numerical Mathematics}
\centerline{Faculty of Mathematics, Physics and Informatics of the Comenius University}
\centerline{Mlynsk\' a dolina, 842 48 Bratislava, Slovakia}

\bigskip
\centerline{$^\spadesuit$Academy for Multidisciplinary studies, Capital Normal University}
\centerline{ West 3rd Ring North Road 105, 100048 Beijing, P.~R.~China}


\begin{abstract}
In this paper we study the convergence rate of a finite volume approximation of the compressible Navier--Stokes--Fourier system.
To this end we first show the local existence of a highly regular unique strong solution and
analyse its global extension in time as far as the density and temperature remain bounded.
We make a physically reasonable assumption that the numerical density and temperature are uniformly bounded from above and below. The relative energy provides us an elegant way to derive
a priori error estimates between finite volume solutions and the strong solution.
\end{abstract}



\section{Introduction}
Numerical simulation of gas dynamics plays an important role in a wide range of modern industrial and real life applications, such as vehicle engineering, engines design, aerospace or weather forecast. In the past decades, lots of effort was put in the development of efficient (high order, high resolution and parallel) and accurate (stable and convergent) numerical methods. Comparing to the great success in the development of reliable numerical methods, the progress in their analysis is much slower.

In this paper, we study the convergence rate of a finite volume (FV) method proposed by Feireisl et al.~\cite{FLMS_FVNSF} for the compressible Navier--Stokes--Fourier system.
Our main tool is the relative energy functional originally proposed by Dafermos~\cite{Dafermos} in order to measure the distance of a weak and a classical solution. The idea has been adapted to study the error between a numerical solution and a classical solution, see Gallou\"et et al. \cite{GallouetMAC} and Feireisl et al. \cite{FHMN}. These results are based on a discrete version of the relative energy functional, where the target smooth solution is projected into the discrete space of numerical solutions. Recently, Feireisl et al. \cite{FLS} showed that using the continuous version of the relative energy functional yields a better convergence rate compared to its discrete version. This new strategy directly employs the consistency error of the numerical scheme and avoids lengthy calculations in the discrete space.
This approach has been already successfully applied for the error estimates of the compressible Euler equations, see Feireisl et al. \cite{FeLMMiSh} for the barotropic case and see Luk\'{a}\v{c}ov\'{a} et al. \cite{LMSY} for the full Euler system.
In the present paper we use the \emph{continuous form of the relative energy} and study {\em the convergence rate of a FV method} for the Navier--Stokes--Fourier system.  To this end, we need a new refined version of the consistency errors as well as stability estimates. As far as we know, this is the first result in the literature, where the error estimates of a fully discrete numerical scheme for the compressible viscous heat-conducting fluids are derived. In this paper we concentrate on theoretical error analysis, numerical simulations illustrating reliability of the studied FV method were shown in \cite{LMS21}.

The Navier--Stokes--Fourier system describes the motion of compressible, viscous, and heat conducting fluids
\begin{subequations}\label{PDE_nsf}
\begin{equation}
\pd_t \vr  + \Div (\vr \vu ) = 0,
\end{equation}
\begin{equation}
\pd_t (\vr \vu)  + \Div (\vr \vu \otimes \vu  ) + \Grad p= \Div \bS ,
\end{equation}
\begin{equation}
\pd_t (\vr e)  +  \Div ( \vr e\vu) - \Div (\kappa \Grad \vt) =  \bS: \Grad \vu  -p \Div \vu
\end{equation}
\end{subequations}
in the time--space cylinder $(0,T)\times \Omega$,  $\Omega\subset \R^d$, $d=2, 3$.
Here, $\vr, \vu, \vt, p, e$ are the fluid density, velocity,  absolute temperature, pressure and internal energy, respectively. Further, $\kappa>0$ is a heat conductivity constant and $\bS$ is the viscous stress tensor
 $$\bS=\bS(\Grad \vu)= 2\mu \bDu  + \lambda \Div \vu \I \quad \mbox{with} \quad  \bDu=\frac{\Grad \vu + \Grad^T \vu}{2},$$
 where $\mu>0$ and $ \lambda \geq 0$ are constant viscosity coefficients.
We consider the standard pressure law of a perfect gas
\[ p = (\gamma-1) \vr e,\ e=c_v \vt,  \]
where $\gamma>1$ is the adiabatic coefficient and $c_v = \frac{1}{\gamma-1}$ is the specific heat per constant volume.
For simplicity, we consider periodic boundary conditions and identify the fluid domain with a flat torus, i.e., $\Omega=\Of=\left([0,1]_{[0,1]}\right)^d$.
To close the system we impose the initial conditions
\begin{equation}\label{INI}
(\vr (0), \vu (0),  \vt (0))=(\vr _0, \vu _0,\vt _0)  \text{ with } \vr_0 >0\text{ and } \vt_0>0.
\end{equation}
In view of the second law of thermodynamics, any solution to the system \eqref{PDE_nsf}--\eqref{INI} satisfies the entropy balance
\begin{equation}\label{pde_s}
\vr( \pd_t s + \vu \cdot \Grad s) - \Div \left(\frac{\kappa \Grad \vt}{\vt} \right)
= \frac{1}{\vt} \left(\bS(\Grad \vu): \Grad \vu + \kappa \frac{|\Grad \vt|^2}{\vt}  \right) \geq 0,
\end{equation}
where $s$ is the physical entropy
\[s(\vr, \vt)= \log \left(\frac{\vt^{c_v}}{\vr} \right)  \qquad \mbox{ for } \vr > 0,\ \vt >0.
\]
The rest of the paper is structured as follows.
In Section \ref{SS} we show the existence of a strong solution.
Next, we introduce a FV method in Section.~\ref{sec_NT}.
Then, we show the stability and consistency of the FV approximation in Section~\ref{sec_CS}. Finally, we present the convergence rate of the FV  method in Section~\ref{sec_EE}.

\section{Strong solution}\label{SS}
In this section, we focus on the existence results to the Navier--Stokes--Fourier system \eqref{PDE_nsf} in the class of strong solutions. In particular, we emphasize that the classical solution inherits the $W^{k,2}$-regularity from the initial data for any $k\geq 4$, as a consequence of the fact that \eqref{PDE_nsf} can be seen as a symmetric hyperbolic-parabolic system, cf. Theorem \ref{local ex}. Moreover, due to the blow-up criterion established by Huang, Li \cite{HuaLi} and by Feireisl, Wen, Zhu \cite{FeiWenZhu}, the solution exists globally in time, as long as it is bounded.

\begin{Proposition}  \label{strong sol}
	Let $k\geq 6$ and let the initial data \eqref{INI} be such that
	\begin{equation*}
	\vr_0>0, \quad (\vr_0, \vu_0, \vt_0) \in W^{k,2}(\Of; \mathbb{R}^{1+d+1}).
	\end{equation*}
	Let $(\vr, \vu, \vt)$ be a weak solution to problem \eqref{PDE_nsf}, \eqref{INI} such that
	\begin{equation} \label{bound vr vt}
	0< \Un{\vr} \leq \vr \leq \Ov{\vr}, \quad 0< \Un{\vt} \leq \vt \leq \Ov{\vt} \quad \mbox{a.e. in} \ (0,T)\times \tor.
	\end{equation}
	Then $(\vr, \vu, \vt)$ is a strong solution of \eqref{PDE_nsf}, \eqref{INI} in $[0,T] \times \Of$ satisfying
	\begin{equation*}
		(\vr, \vu, \vt) \in C([0,T]; W^{k,2}(\Of; \mathbb{R}^{1+d+1})).
	\end{equation*}
	Moreover, the following estimate
	\begin{equation}\label{STC}
	\begin{aligned}
	\| (\vr, \vu, \vt) \|_{C([0,T]; W^{k,2}(\Of))} + \| (\vr, \vu, \vt) \|_{C^1([0,T] \times \Of)}& \\[0.1cm]
	+ \| (\pd_{t}\vr, \pd_{t}\vu, \pd_{t}\vt) \|_{C([0,T] \times \Of)} + \| (\pd_{t}^2\vr, \pd_{t}^2\vu, \pd_{t}^2\vt) \|_{C([0,T] \times \Of)} &\leq C \| (\vr_0, \vu_0, \vt_0)\|_{W^{k,2}}
	\end{aligned}
	\end{equation}
	holds, where the positive constant $C$ depends solely on $T, \Un{\vr}, \Ov{\vr}, \Un{\vt}, \Ov{\vt}$ and $\| (\vr_0, \vu_0, \vt_0)\|_{W^{k,2}}$.
\end{Proposition}

\begin{proof}
	Applying the local-existence result developed in Theorem \ref{local ex}, there exists $T_1$, $0 < T_1 \leq T$, such that a unique strong solution $(\tvr, \tvu, \tvt)$ to problem \eqref{PDE_nsf}, \eqref{INI} exists on $[0,T_1]$ and
	\begin{equation*}
	(\tvr, \tvu, \tvt) \in C([0,T_1]; W^{k,2}(\Of; \mathbb{R}^{1+d+1})).
	\end{equation*}
	Moreover, from the Sobolev embedding $W^{k-4,2}(\Of) \hookrightarrow C(\Of)$ whenever $k\geq 6$, from \eqref{estimate ss} we can deduce that there exists a positive constant $C=C(T_1)$ such that
	\begin{equation} \label{first estimate}
	\begin{aligned}
	\| (\tvr, \tvu, \tvt) \|_{C([0,T_1]; W^{k,2}(\Of))} &+ \| (\tvr, \tvu, \tvt) \|_{C^1([0,T_1] \times \Of)} \\[0.1cm]
	+ \| (\pd_{t}\tvr, \pd_{t}\tvu, \pd_{t}\tvt) \|_{C([0,T_1] \times \Of)} &+ \| (\pd_{t}^2\tvr, \pd_{t}^2\tvu, \pd_{t}^2\tvt) \|_{C([0,T_1] \times \Of)} \leq C \| (\vr_0, \vu_0, \vt_0)\|_{W^{k,2}}.
	\end{aligned}
	\end{equation}
	Now, from the weak-strong uniqueness principle proven in \cite[Theorem 2.1]{FeiNov10}, the weak solution $(\vr, \vu,\vt)$ must coincide with the strong solution $(\tvr, \tvu, \tvt)$ in $[0,T_1] \times \Of$. According to the blow-up criterion developed by Feireisl, Wen and Zhu \cite[Theorem 1.4]{FeiWenZhu}, if $T_1 < T$ were the maximal time of existence, then
	\begin{equation} \label{blow up}
	\limsup_{t \rightarrow T_1^-} \left( \| \vr(t)\|_{L^{\infty}}+ \| \vt(t)\|_{L^{\infty}}  \right) = \infty,
	\end{equation}
	which would lead to a contradiction due to hypothesis \eqref{bound vr vt}. 
	Consequently $T_1 = T$ which concludes the proof.
\end{proof}

\begin{Remark}
	We point out that the compatibility condition (1.10) required for the initial data in \cite{FeiWenZhu} is not necessary in this context as we are not allowing the presence of vacuum zones.
\end{Remark}
\begin{Remark}
	Proposition \ref{strong sol} still holds if we replace the boundedness of temperature with the one of velocity; specifically, instead of having $\vt \leq \Ov{\vt}$ in \eqref{bound vr vt}, we could require $|\vu| \leq \Ov{\vu}$. Then, instead of \eqref{blow up}, we could use the Serrin-type blow-up criterion developed by Huang and Li \cite[Theorem 1.2]{HuaLi}: if $T_1$ is the maximal time of existence, then
	\begin{equation*}
	\limsup_{t \rightarrow T_1^-} \left( \| \vr(t)\|_{L^{\infty}}+ \| \vu(t)\|_{L^{\infty}}  \right) = \infty.
	\end{equation*}
\end{Remark}

\section{The finite volume method}\label{sec_NT}
In this section we introduce the finite volume method originally proposed in \cite{FLMS_FVNSF}.
\subsection{Notations}
\paragraph{Mesh.}  Let $\grid$ be a uniform structured mesh with the mesh size $h$ formed by squares in 2D, and cubes in 3D such that $\tor = \bigcup_{K \in \grid} K$.

We denote by $\faces$ the set of all faces of $\grid$, and by $\facesi$ the set of all faces that are orthogonal to the $\ith$ basis vector $\ve_i$ of the canonical system.  The set of all faces of an element $K$ is denoted by $\facesK$.
Moreover, we write $\sigma= K|L$ if $\sigma \in \faces$ is the common face of elements $K$ and $L$. Further, we denote by $|K|$  and $|\sigma|$ the Lebesgue measure
of an element $K\in\grid$ and a face $\sigma \in \faces$, respectively.
Moreover, $x_K$ and $x_\sigma$ stand for the barycenter of $K$ and $\sigma$, respectively.

In what follows we shall use an abbreviated notation for the integrals,
\begin{align*}
\intTor{\bdot} := \sum_{K\in\grid}\intK{\bdot} \quad \mbox{and} \quad \intfacesint{\bdot}:=\sum_{\sigma\in\facesint}\intSh{\bdot}.
\end{align*}

\begin{figure}[hbt]
\centering
\begin{tikzpicture}[scale=1.0]
\draw[-,very thick](0,-2)--(4,-2)--(4,2)--(0,2)--(0,-2)--(-4,-2)--(-4,2)--(0,2);
\draw[blue!90,fill=blue!20,very thick] (0,-2)--(2.0,-2)--(2.0,2)--(0,2)--(0,-2);
\draw[-,very thick, red=90!, pattern=north west  lines, pattern color=red!30] (0,-2)--(0,2)--(-2.0,2)--(-2.0,-2)--(0,-2);

\path node at (-3.0,0) { $K$};
\path node at (3.0,0)  { $L$};
\path node at (-2.0,0) {$ \bullet$};
\path node at (-2.3,-0.3) {$ x_K$};
\path node at (2.0,0) {$\bullet$};
\path node at (2.2,-0.3) {$ x_L$};
\path node at (0,0) {$\bullet$};
\path node at (0.3,-0.3) {$ x_\sigma$};

\path (-0.4,0.8) node[rotate=90] { $\sigma=\overrightarrow{K|L}$};
 \path (-1.5,1.4) node[] { $D_{\sigma K}$};
 \path (1.5,1.4) node[] { $D_{\sigma L}$};
\end{tikzpicture}
\caption{Dual cell $D_{\sigma} = D_{\sigma K} \cup D_{\sigma L}$ for 2D structured mesh.}\label{figmesh}
\end{figure}
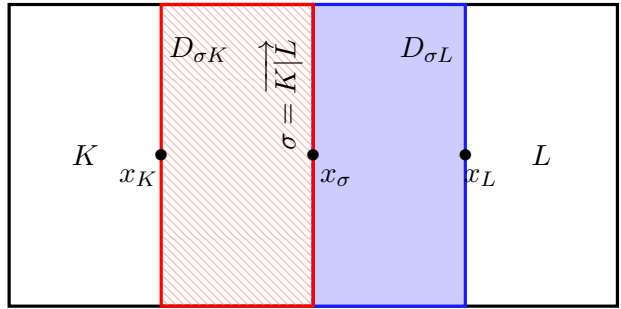

\paragraph{Dual grid.} For any $\sigma=K|L\in \facesi$, we define a dual  cell  $D_\sigma := D_{\sigma,K} \cup D_{\sigma,L}$, where the  polygon $D_{\sigma,K}$ (resp.~$D_{\sigma,L}$) is  a half of $K$ (resp. $L$), i.e.,
\begin{align*}
& D_{\sigma,K} = \left\{ x\in K\mid  x_i \in \co{(x_K)_i}{(x_\sigma)_i} \right\},
\\
& \co{A}{B} \equiv [ \min\{A,B\} , \max\{A,B\}],
\end{align*}
see Figure~\ref{figmesh} for a two-dimensional example of dual cells. Further, we denote the $\ith$ dual grid $\Di$ as
$$
\Di =\bigcup_{\sigma \in \facesi} D_{\sigma}.
$$

 \paragraph{Function space.}
 The symbol $Q_h$ stands for the set of piecewise constant functions on the grid $\grid$.
Note that hereinafter $\vv_h\in Q_h$ means that every component of a vector--valued function $\vv_h$ belongs to the set $Q_h$. We define the following  projection operators
\begin{equation*}
\Piq  \phi (x) = \sum_{K \in \grid}  \frac{1_{K}(x)}{|K|} \int_K \phi \dx , \quad
\Piwi \phi(x) = \sum_{\sigma \in \facesi} \frac{1_{D_\sigma}(x)}{|\sigma|} \int_{\sigma} \phi \ds,
\end{equation*}
where
$1_{D}(x) = \begin{cases}1 & x \in D ,\\0 & \mbox{otherwise,}\end{cases}$ for $D=K, D_\sigma$ being an element of $\grid$ or $\Di$, respectively.
Note that for any $\phi \in W^{1,\infty}(\tor)$
\begin{equation}\label{proj}
\norm{\phi - \Piwi \phi}_{L^\infty(\tor)} \leq h \norm{\phi}_{W^{1,\infty}(\tor)}.
\end{equation}
The average and jump operators are denoted as
\[
\avs{v}(x) = \frac{v^{\rm in}(x) + v^{\rm out}(x) }{2},\quad
\jump{ v }  = v^{\rm out}(x) - v^{\rm in}(x),
\]
 where
 \[
v^{\rm out}(x) = \lim_{\delta \to 0+} v(x + \delta \vc{n}),\quad
v^{\rm in}(x) = \lim_{\delta \to 0+} v(x - \delta \vc{n}),\
\]
whenever $x \in \sigma \in \facesint$  and $\vc{n}$ is the outer normal vector to $\sigma$.

\noindent {\bf Discrete operators.} The discrete gradient $\Gradh$ and divergence $\Divh$ operators acting on the mesh $\grid$ are defined as  follows,
\begin{equation*}
\begin{aligned}
\left( \Gradh r\right)_K  &=  \frac{|\sigma|}{|K|} \sum_{\sigma \in \facesK} \avs{r} \vc{n}, & \left( \Divh \bfv\right)_K &=
 \frac{|\sigma|}{|K|} \sum_{\sigma \in \facesK} \avs{\bfv} \cdot \vc{n} ,
\\
 \Gradh r(x)  &=  \sum_{K \in \grid} \left( \Gradh r\right)_K  1_K{(x)}, & \Divh \bfv(x) &= \sum_{K\in\grid} \left( \Divh \bfv\right)_K 1_K{(x)} .
\end{aligned}
\end{equation*}
We point out that $\Gradh$ and $\Divh$ are  discrete variants of $\Grad$ and $\Div$ on the cell $K$, respectively. Moreover, we define
\begin{align*}
\Gradh \bfv = \left( \Gradh v_{1}, \dots,   \Gradh v_{d}\right)^T, \quad
\bD_h (\bfv) = \frac{\Gradh \bfv   +\Gradh^T \bfv}{2}
\quad \mbox{yielding} \quad
\Divh \bfv = \mbox{tr} (\Gradh \bfv) =   \mbox{tr} (\bD_h (\bfv))
\end{align*}
and the following difference operator 
for any $r_h \in Q_h$
\begin{equation}\label{grad_edge}
\begin{aligned}
&  \Gradedge r_h(x)  = \sum_{\sigma\in\faces}  \left(\Gradedge r_h \right)_{\sigma}1_{D_\sigma}{(x)},  \quad \left(\Gradedge r_h\right) _{\sigma} =\frac{\jump{r_h} }{ h } \vc{n}.
\\&
\Delta_h r_h(x) =  \sum_{K \in \grid} \left(\Delta_h r_h\right)_K  1_K{(x)}, \quad  \left(\Delta_h r_h\right)_K =  (\Divh \Gradedge r_h)_K = \frac{|\sigma|}{|K|} \sum_{\sigma \in \facesK} \frac{ \jump{r_h} }{h}.
\end{aligned}
\end{equation}

It is crucial that the following duality relationships between discrete operators hold
\begin{subequations}\label{dis_op}
\begin{align}
& \intO{ r_h\Divh \vv_h } =  -\intO{ \Gradh r_h \cdot \vv_h },
\\
&  \intO{ \Delta_h r_h \cdot f_h } =  -\intO{ \Gradedge r_h \cdot  \Gradedge f_h } =  \intO{r_h \cdot  \Delta_h  f_h },  \label{dis_op3}
\\
&\intO{ r_h\Div \bfphi } = - \intO{ \Gradedge r_h \cdot  \Piw \bfphi } \quad \mbox{for} \  \bfphi \in W^{1,1}(\tor;\R^d),
\end{align}
\end{subequations}
where $\vv_h, r_h, f_h \in Q_h$ and $\Piw \bfphi=(\Piw^{(1)}\phi_1, \ldots, \Piw^{(d)}\phi_d)$. We refer to \cite[Lemma 5 \& Lemma 7]{FeLMMiSh} for proofs.

\paragraph{Time discretization.}
Given a time step $\TS >0$ we divide the time interval $[0,T]$ into $N_T=T/\TS$ uniform parts  and denote $t^k= k\TS.$ Then $v_h^k \approx v_h(t^k)$ is the approximation of a function $v_h$ at time $t^k,$  $k=0,1,\ldots,N_T$.
By $v_h \in L_{\TS}(0,T; Q_h)$ we denote a piecewise constant  in time function
\begin{equation*}
v_h(t, \cdot) =v_h^0  \ \mbox{ for } \  t < \Delta t, \quad  v_h(t,\cdot)=v_h^k \ \mbox{ for } \ t\in [k\TS,(k+1)\TS), \ \   k=1,\cdots, N_T.
\end{equation*}
Further, we define the discrete time derivative by
\[
 D_t v_h(t) = \frac{v_h (t) - v_h^\triangleleft}{\TS} \quad \mbox{with} \quad v_h^\triangleleft  = v_h (t - \TS).
\]

\subsection{Numerical scheme}
We  are now ready to propose a FV method for the compressible Navier--Stokes-Fourier system \eqref{PDE_nsf} in a weak form, see also \cite{FLMS_FVNSF}.

\begin{Definition}[FV method]
  Given the initial data  \eqref{INI} we set $(\vrh^0,\vuh^0, \vth^0) =(\Piq\vr_0, \Piq\bm{u}_0, \Piq \vt_0)$. We say that the triple
$(\vrh,\vuh, \vth) \in L_{\TS}(0,T; Q_h \times \bQh \times Q_h)$ is a finite volume approximation of the Navier--Stokes--Fourier system \eqref{PDE_nsf}, \eqref{INI}, if the following system of algebraic equations holds
\begin{subequations}\label{scheme}
\begin{align}\label{scheme_D}
&\intO{ D_t \vrh  \phi_h } - \intfacesint{  \Fup (\vrh ,\vuh ) \jump{\phi_h}   } = 0 \quad \mbox{for all}\ \phi_h \in Q_h,
\\  \label{scheme_M}
&\intO{ D_t  (\vrh  \vuh ) \cdot \bfphi_h } - \intfacesint{ \Fup  (\vrh  \vuh ,\vuh ) \cdot \jump{\bfphi_h}   }
+\intO{ (\bS_h - p_h \I):\Gradh \bfphi_h } =0
\quad \mbox{for all } \bfphi_h \in \bQh,
\\ \label{scheme_T}
&c_v\intO{ D_t (\vrh  \vth ) \phi_h } - c_v\intfacesint{  \Fup (\vrh \vth ,\vuh )\jump{\phi_h} }
+\intfacesint{  \frac{\kappa}{ h } \jump{\vth}  \jump{ \phi_h}  }
\br
 & \hspace{7cm}= \intO{  (\bS_h-p_h \I ):\Gradh \vuh  \phi_h}, \quad \mbox{for all}\ \phi_h \in Q_h ,
\end{align}
\end{subequations}
where $\bS_h = 2 \mu \Dhuh  + \lambda  \Divh   \vuh \I$,  the discrete pressure $p_h$ is given by
\begin{equation*}
p_h(t)= \vrh(t) \vth(t)  \quad \mbox{ for } \vth(t) >0  \quad  \mbox{ and } \quad p_h(t) = 0 \quad \mbox{ for } \vth(t) \leq 0,
\end{equation*}
and $\Fup$, $\eps \in (-1,1)$, is the so-called diffusive upwind flux
\begin{align*}
\Fup (r_h,\vuh)
=\Up[r_h, \vuh] - \muh \jump{ r_h }, \;
\Up [r_h, \vuh]
= \avs{r_h} \ \auh \cdot \vn - \frac{1}{2} |\auh \cdot \vn| \jump{r_h}.
\end{align*}
\end{Definition}

\begin{Remark}
Note that the existence of a solution $(\vrh, \vuh, \vth)$ to the nonlinear algebraic system \eqref{scheme} can be done using the fixed point theorem \cite[Theorem 15]{FeLMMiSh} in the same way as in \cite[Lemma 11.3]{FeLMMiSh}.
In our recent paper \cite{LMS21} we have presented several numerical experiments obtained by the FV method \eqref{scheme} that illustrate its reliability and robustness.
\end{Remark}

\section{Stability and consistency of the FV method}\label{sec_CS}
In this section we discuss the stability and consistency of the FV method
\eqref{scheme}. To begin, we recall from \cite{FLMS_FVNSF} the following structure-preserving properties of the method.
\begin{Lemma}\label{lm_pro}
Let $(\vrh,\vuh,\vth)$ be a numerical solution of the FV method \eqref{scheme} with positive initial data  $\vrh^0>0$ and $\vth^0>0.$
Then we have
\begin{itemize}

\item {\bf \textit {Positivity of the density.}} 
\begin{equation*}
\vrh(t) > 0 \quad  \mbox{for all} \ \ t\in(0,T).
\end{equation*}

\item {\bf \textit{ Positivity of the temperature.}} 
\begin{equation*}
\vth(t) > 0 \quad  \mbox{for all} \ \ t\in(0,T).
\end{equation*}

\item {\bf \textit{Mass conservation.}}
	\begin{equation*}
	\intO{\vrh(t)} = \intO{\vrh(0)} = \intO{\vr_0}>0.
	\end{equation*}
\item {\bf \textit{Energy balance.}}
	\begin{align}\label{energy_stability}
	 &D_t \intO{ \left(\frac{1}{2}  \vrh  |\vuh |^2 +  c_v \vrh  \vth \right) }
	 +  h^\eps \intfacesint{  \avs{ \vrh  }  \abs{\jump{\vuh}}^2 }
	\br
	&\hspace{3cm}+ \frac{\TS}{2} \intO{ \vrh^\triangleleft|D_t \vuh |^2  }
	+ \frac12 \intfacesint{ \vrh^{\rm up} |\avs{\vuh } \cdot \vc{n} | \abs{\jump{\vuh}}^2   }
	= 0.
	\end{align}
\item {\bf \textit{Entropy balance.}}
	\begin{align}\label{eq_entropy_stability}
	 \intO{ D_t \left(\vrh  s_h  \right) \phi_h  } & -
	\intfacesint{ \Up(\vrh  s_h , \vuh ) \jump{\phi_h} }
	+ \intfacesint{\frac{\kappa}{ h } \jump{\vth}   \jump{ \frac{\phi_h}{\vth }} }
	\br
	& - \intO{ \bS_h:\Gradh \vuh \frac{\phi_h}{\vth } }
	 =  D_s(\phi_h)  + R_{s}(\phi_h),
	\end{align}
	where $\phi_h \in Q_h$, $D_s(\phi_h) = D_1 (\phi_h)  + D_2(\phi_h)+ D_3(\phi_h)$  and
\begin{equation}\label{entropy_dissipation}
\begin{aligned}
D_{1}(\phi_h)  := & \  \TS \intO{ \left( \frac{ | D_t \vr _h |^2 }{2 \xi _{\vr,h} }  +\frac{c_v \vrh^\triangleleft  |D_t \vth |^2 }{2 |\xi_{\vt,h} |^2 } \right)  \phi_h},
\\
 D_{2}(\phi_h)  := & \ \frac12 \intfacesint{    \phi_h^{\rm down}  \abs{ \avs{\vuh}\cdot \vc{n} }      \jump{ (\vrh, p_h) } \cdot \nabla_{(\vr,p)}^2 (-\vr s)|_{\vv_1^*}\cdot \jump{ (\vrh, p_h) } },
\\
D_{3}(\phi_h) :=
 & \ h^{\eps} \intfacesint{     \avs{\phi_h}  \jump{ (\vrh, p_h) } \cdot \nabla_{(\vr,p)}^2 (-\vr s)|_{\vv_2^*}\cdot \jump{ (\vrh, p_h) }},
 \\
R_{s}(\phi_h) := & \  h^{\eps}  \intfacesint{    \jump{\phi_h}  \cdot  \bigg( \avs{  \nabla_\vr(-\vrh s_h) }  \jump{ \vrh}
	+ \avs{\nabla_p(-\vrh s_h)} \jump{ p_h}\bigg)},
\\
  \xi_{\vr,h}  \in &\co{\vrh^{\triangleleft}}{\vrh}, \
\xi_{\vt,h} \in \co{\vth^{\triangleleft}}{\vth}, \
\vv_1^*, \vv_2^* \in \co{(\vrh^{\rm in}, p_h^{\rm in})}{(\vrh^{\rm out}, p_h^{\rm out})}.
\end{aligned}
\end{equation}
Moreover, it holds that $D_s(\phi_h) \geq 0$  for $ \phi_h \geq 0$.
\end{itemize}
\end{Lemma}

\begin{Remark}
We point out that the entropy Hessian matrix $\nabla_{(\vr,p)}^2 (-\vr s)$  is symmetric positive definite for $(\vr,\vt)\in(0,\infty)^2$, i.e.
\begin{equation}\label{EH}
	\nabla_{(\vr,p)}^2 (-\vr s) =
	\begin{pmatrix}
	\frac{1+c_v}{\vr}  & -\frac{c_v}{\vr \vt} \\
	 -\frac{c_v}{\vr \vt} & \frac{c_v}{\vr \vt^2}
	\end{pmatrix} >0.
\end{equation}
Assuming that $\vrh$ and $\vth$ are uniformly bounded from above and below by positive constants, see \eqref{AS} below, we obtain that the entropy Hessian matrix $\nabla_{(\vr,p)}^2 (-\vr s)(\vrh, \vth) $ is uniformly bounded from above and below by some positive constants, see Appendix \ref{sec-EH} for details.
\end{Remark}

With the boundedness of the entropy Hessian matrix $\nabla_{(\vr,p)}^2 (-\vr s)(\vrh, \vth)$, we obtain the following uniform bounds from the entropy balance \eqref{eq_entropy_stability} with test function $\phi_h = 1$.

\begin{Corollary}[Uniform bounds]
Let $(\vrh ,\vuh ,\vth )$ be a  numerical solution of the FV method \eqref{scheme} with positive initial data  $\vrh^0>0$ and $\vth^0>0.$
	Under the assumption
	\begin{equation}\label{AS}
		  0< \Uvt \leq \vth \leq \Ovt ,\ 0< \underline{\vr} \leq \vrh \leq \Ov{\vr} \ \mbox{ uniformly for }\ h \rightarrow 0,
	\end{equation}
	there exist some constants $\underline{s} \leq \Ov{s}$, $0<\underline{p} \leq \Ov{p}$ such that the entropy and pressure are bounded,
		\begin{equation}\label{ESen}
		\underline{s} \leq  s_h \leq \Ov{s}  \mbox{ and }  \underline{p} \leq p_h \leq \Ov{p}  \ \mbox{ uniformly for }\  h \rightarrow 0.
	\end{equation}
Moreover, we have the following estimates
\begin{subequations}\label{N}
\begin{equation}\label{N1}
\norm{\vuh }_{L^\infty (0,T; L^{2}(\tor)) } +  \norm{\Gradedge \vth}_{L^2((0,T)\times\tor)}  + \norm{\bS_h}_{L^2((0,T)\times\tor)}  \aleq 1 ,
\end{equation}
\begin{equation}\label{N2}
\int_0^{\tau}\intfacesint{   \left( h^{\eps} + \abs{ \avs{\vuh}\cdot \vc{n} } \right) \left( \jump{\vrh}^2 + \jump{p_h}^2 + \abs{\jump{\vuh}}^2\right)      }\dt \aleq 1.
\end{equation}
\end{subequations}
\end{Corollary}
\begin{Remark}
It is easy to verify that
\begin{equation}\label{N3}
h^{\eps+1} \norm{\Gradedge \vuh}_{L^2L^2}^2 = h^\eps\intfacesint{   \abs{\jump{ \vuh }} ^2   }  \aleq 1
\end{equation}
and $\displaystyle{ \intTO{\Gradh \vuh : \Gradh^T \vuh } = \norm{\Divh \vuh}_{L^2L^2}^2 \geq 0}$, which further implies
\begin{equation}\label{N4}
\norm{\Gradh \vuh}_{L^2L^2}   \approx \norm{\Dhuh}_{L^2L^2}  \approx \norm{\bS_h}_{L^2L^2}  \aleq 1.
\end{equation}
Note that the notation ``$a \approx b$" means $a \aleq b$ and $b \aleq a$.
\end{Remark}

Now we are ready to analyse the consistency of the FV method \eqref{scheme}.
\begin{Lemma}[Consistency formulation] \label{lem_C1}
Let $(\vrh, \vuh, \vth)$ be a numerical solution obtained by the FV scheme \eqref{scheme}
with  $(\TS,h) \in (0,1)^2$, $-1 < \eps <1$.
Let the assumption \eqref{AS} hold. Then,
\begin{subequations}\label{cP}
\begin{equation} \label{cP1}
\left[ \intO{ \vrh\phi }   \right]_{t=0}^\tau=
 \intTOB{  \vrh \partial_t \phi + \vrh \vuh \cdot \Grad \phi }   +   e_\vr(\phi,\TS, h,\tau)
\end{equation}
for any $\phi \in L^\infty(0,T;W^{2,\infty}(\tor)),$ $\partial_t\phi, \partial^2_t\phi \in L^\infty((0,T) \times  \tor)$  ;
\begin{align} \label{cP2}
\left[ \intO{ \vrh \vuh \cdot \bfphi } \right]_{t=0}^\tau  & =
 \intTOB{  \vrh \vuh \cdot \partial_t \bfphi + \vrh \vuh \otimes \vuh  : \Grad \bfphi}
\br
&+ \intTO{(p_h\I - \bS_h): \Grad \bfphi}
+  e_{\vm}(\bfphi, \TS, h,\tau)
 \end{align}
for any $\bfphi \in L^\infty(0,T;W^{2,\infty}(\tor;\mathbb R^d)),$ $\partial_t\bfphi, \partial^2_t\bfphi \in L^\infty((0,T) \times  \tor; \mathbb R^d)$;
\begin{equation}\label{cP3}
\begin{aligned}
&\left[ \intO{ \vrh s_h \phi } \right]_{t=0}^\tau  \geq
 \intTO{  \vrh s_h (\pd_t\phi + \vuh \cdot \Grad \phi ) } +  \intTO{ \bS_h:\Gradh \vuh \frac{ \phi}{\vth} }
\\
&\hspace{2cm}+   \intTO{ \frac{ \kappa \phi}{ \vthout \vth} \abs{\Gradedge \vth }^2 }
-   \intTO{ \frac{\kappa}{\vth} \Gradedge \vth \cdot \Grad \phi }
+  e_{s}(\phi, \TS, h,\tau)
\end{aligned}
\end{equation}
for any  $\phi \in L^\infty(0,T;W^{2,\infty}(\tor)),$ $\partial_t\phi, \partial^2_t\phi \in L^\infty((0,T) \times  \tor)$, $\phi \geq 0$.
\end{subequations}
Here
\begin{subequations}\label{ecs}
\begin{align}\label{ec1}
 \Abs{e_\vr(\phi,\TS, h,\tau) } &  \leq C_\varrho\left(   \TS +h+h^{1-\eps} + h^{1+\eps}\right),
\\
 \Abs{ e_{\vm}(\bfphi, \TS, h,\tau) } & \leq  C_{\vm}\left(\TS+h+h^{(1-\eps)/2}+h^{1+\eps}\right), \label{ec2}
\\
 \Abs{ e_{s}(\phi, \TS, h,\tau) } & \leq  C_s\left( \TS +h+h^{1-\eps} + h^{(1+\eps)/2}\right)\label{ec3}
\end{align}
and the constants $C_\vr,$ $C_{\vm}$, $C_s$ are independent of the discretization parameters $h,$ $\TS.$
\end{subequations}
\end{Lemma}
\noindent The proof of the above lemma has been done in \cite{FLMS_FVNSF} with more regular test functions but without a clear order of the consistency error.
Here we present an optimal consistency errors with minimum regularity assumptions on the test functions. We note that in the next section the exact strong solution will be used as a test function for the error estimates.
For completeness we provide the details of the proof in Appendix~\ref{app_cs}.

\section{Error estimates}\label{sec_EE}
In this section we show the main result, that is the error estimate of the FV approximation \eqref{scheme}. To begin, we recall the relative energy functional which measures the distance between the numerical solution $(\vrh,\vuh,\vth)$ and the strong solution $(\tvr,\tvu,\tvt)$
\begin{align*}
\RE{\vrh,\vuh,\vth}{\tvr,\tvu,\tvt}
= \intOB{\frac12 \vrh |\vuh - \tvu|^2 + \EH{\vrh,\vth}{\tvr,\tvt}}
\end{align*}
with
\begin{align*}
\EH{\vr,\vt}{\tvr,\tvt} = \Hc(\vr,\vt)  - \frac{\pd \Hc(\tvr,\tvt)}{\pd \vr}(\vr -\tvr) - \Hc(\tvr,\tvt) \quad \mbox{and} \quad
\Hc(\vr,\vt) = \vr( c_v   \vt  -   \vT s(\vr,\vt)   ).
\end{align*}
Then it is straightforward to check
\begin{equation}\label{HP}
\begin{aligned}
&\frac{\pd \Hc(\vr,\vt)}{\pd \vr} =
  c_v \vt - \vT \big( s(\vr,\vt)  -1 \big),
\quad
\tvr \frac{\pd \Hc(\tvr,\tvt)}{\pd \vr}  -   \Hc(\tvr,\tvt)
= \tp ,
\\&
 \pd_y \frac{\pd \Hc(\tvr,\tvt)}{\pd \vr} + \ts \pd_y \tvt =  \frac{1}{\tvr} \pd_y \tp, \quad y=t,x,
\end{aligned}
\end{equation}
where we have denoted $\tp = p(\tvr, \tvt) ,\  \ts =s(\tvr,\tvt)$.

\begin{Remark}
Let $\eta = \vr s$. Taking $(\vr, \eta)$ as independent variables, we reformulate $\EH{\vr,\vt}{\tvr,\tvt}$ as
\begin{align*}
\EH{\vr,\vt}{\tvr,\tvt}= \vr e  - \frac{\pd  (\vr e)}{\pd \vr}\Big|_{(\tvr, \widetilde{\eta})} (\vr - \tvr) -  \frac{\pd  (\vr e)}{\pd \eta}\Big|_{(\tvr, \widetilde{\eta})}  (\eta - \widetilde{\eta} ) -  \tvr \tilde{e}.
\end{align*}
Due to the convexity of $\vr e$ with respect to $(\vr, \eta)$, we know that $\EH{\vr,\vt}{\tvr,\tvt} \geq 0$.
Moreover, it holds under the assumption \eqref{AS} that
\begin{equation}\label{EN}
\begin{aligned}
& \RE{\vrh,\vuh,\vth}{\tvr,\tvu,\tvt} \approx
   \norm{\vuh -\tvu}_{L^2(\tor)}^2 + \norm{ \eta_h  -  \widetilde{\eta} }_{L^2(\tor)}^2 + \norm{\vrh - \tvr}_{L^2(\tor)}^2
\\ &
\approx    \norm{\vuh -\tvu}_{L^2(\tor)}^2 + \norm{ \vth -  \tvt}_{L^2(\tor)}^2 + \norm{\vrh - \tvr}_{L^2(\tor)}^2 +  \norm{s_h - \ts}_{L^2(\tor)}^2.
\end{aligned}
\end{equation}
For more details, we refer to \cite{LMSY} and \cite[Chapter 3.5]{FeiNov_book2}.
\end{Remark}

\begin{Theorem}[Main result: Error estimates]\label{thm_EEB}
Let $\gamma > 1$ and the initial data $(\vr_0, \allowbreak \vu_0,\vt_0)$ satisfy
$$
	\vr_0 >0, \quad (\vr_0, \vu_0, \vt_0) \in W^{k,2}(\Of; \mathbb{R}^{1+d+1}), \quad k\geq 6.
$$
Further, let $(\vrh, \vuh,\vth)$ be a numerical solution obtained by the FV scheme \eqref{scheme}  emanating from the same initial data with $(\TS, h) \in (0,1)^2$ and $-1 < \eps <1$.
In addition, let the numerical density $\vrh$ and temperature $\vth$ be uniformly bounded, i.e. there exist $\Uvt, \Ovt, \Un{\vr}, \Ov{\vr} $ such that
\begin{equation}
		  0< \Uvt \leq \vth \leq \Ovt ,\ 0< \Un{\vr} \leq \vrh \leq \Ov{\vr} \ \mbox{ uniformly for }\ h \rightarrow 0.
\end{equation}

Then there exists a positive number
$$ c=c \left(   T,\Un{\vr},  \Ov{\vr},  \Un{\vt}, \Ov{\vt}, \norm{( \vr_0, \vu_0, \vt_0)}_{W^{k,2}}\right)
$$
such that
\begin{equation}\label{CR1}
\begin{aligned}
&\sup_{0\leq t \leq \tau}\RE{\vrh,\vuh,\vth}{\tvr,\tvu,\tvt}
+\intTO{\frac{\kappa\tvt}{\vth \vthout }|\Gradedge\vth-\Grad\tvt|^2}
\\&\quad +\int_0^\tau \intOB{ \frac{\tvt}{ \vth} \left( 2\mu \left|\Dhuh -  \frac{\vth}{ \tvt} \bD(\tvu) \right|^2  + \lambda  \left| \Divh \vuh -  \frac{\vth}{ \tvt}\Div \tvu \right|^2 \right) }
\\& \leq c( \TS + h^A) \quad \quad \mbox{for all} \ \ \tau \in[0,T],
\end{aligned}
\end{equation}
where $A=  \min\{(1-\eps)/2, (1+\eps)/2\}$ and $(\tvr,\tvu,\tvt)$ is a strong solution satisfying \eqref{STC}.
In particular, we have
\begin{multline}\label{CR2}
\norm{ \vrh - \tvr }_{L^\infty L^2} + \norm{  \vu_h -  \tvu }_{L^\infty L^2}  + \norm{\vth-\tvt}_{L^\infty L^2}
\\ + \norm{ \Gradh \vu_h -  \Grad \tvu }_{L^2 L^2}  + \norm{\Gradedge \vth- \Grad \tvt}_{L^2 L^2}
\leq c( \TS^{1/2} + h^{A/2}).
\end{multline}
\end{Theorem}
\begin{Remark}\hspace{1pt}
\begin{itemize}
\item Recalling Proposition \ref{strong sol}, the choice of initial data with $k\geq 6$ and the assumption \eqref{AS} used in Theorem \ref{thm_EEB} ensure the global existence of strong solution satisfying \eqref{STC}.
\item By choosing $\eps=0$ we obtain $A=1/2$, i.e. the optimal convergence rate in space.
\end{itemize}
\end{Remark}

\begin{proof}[Proof of Theorem \ref{thm_EEB}]
 First, by choosing suitable test functions in the consistency formulation we get
\[
\int_0^\tau  \eqref{energy_stability} \dt
- \eqref{cP3}|_{\phi=\tvt}
+ \eqref{cP1}|_{\phi = \frac{\abs{\tvu}^2}{2} - \frac{\pd \Hc(\tvr,\tvt)}{\pd \vr} }
-\eqref{cP2}|_{\phi = \tvu}
+ \intTO{\pd_t \tp }
\]
and derive the relative energy inequality (see Appendix \ref{app_REI} for the detailed proof)
\begin{equation}\label{REI}
\begin{aligned}
&\left[ \RE{\vrh,\vuh,\vth}{\tvr,\tvu,\tvt}  \right]_0^\tau
 +\intTO{ \frac{\tvt}{ \vth} \left( 2\mu \left|\Dhuh -  \frac{\vth}{ \tvt} \bD(\tvu) \right|^2  + \lambda  \left| \Divh \vuh -  \frac{\vth}{ \tvt}\Div \tvu \right|^2 \right) }
\\
&\hspace{2cm} +\intTO{\frac{\kappa\tvt}{\vth \vthout }|\Gradedge\vth-\Grad\tvt|^2} \leq
R_C   +\sum_{i=1}^5 R_i - R_{\bS} - R_\vt,
\end{aligned}
\end{equation}
where
\begin{align*}
&R_C = e_{\vr}\left( \frac12\abs{\tvu}^2 - \frac{\pd \Hc(\tvr,\tvt)}{\pd \vr} , \TS, h,\tau\right) +e_{\vu} (\tvu, \TS, h,\tau) - e_s(\tvt, \TS, h,\tau),
\\& R_1=
 - \intTO{  \vrh (s_h - \ts)  (\vuh-\tvu) \cdot \Grad \tvt  }  , \quad
R_2= - \int_0^\tau \intO{  \vrh (\vuh - \tvu) \otimes  (\vuh - \tvu): \Grad \tvu },
\\& R_3= \int_0^\tau \intO{   \frac{ \vrh-\tvr}{\tvr} ( \tvu - \vuh )\cdot  \Div \tbS  }
,  \quad
R_4= \intTO{ ( \tp-p_h   -  \pd_\vr \tp (\tvr - \vrh)   -  \pd_\vt \tp (\tvt - \vth)   )\Div \tvu }  ,
\\& R_5=
-\intTO{  \bigg(  (\vrh-\tvr) (s_h- \ts)  + \tvr \big(s_h- \ts - \pd_\vr \ts (\vrh -\tvr)  - \pd_\vt \ts (\vth -\tvt) \big)  \bigg)(\pd_t\tvt + \tvu \cdot \Grad \tvt ) }   ,
\\&
 R_{\bS}  =
  \intTOB{  \vuh \cdot  \Div \tbS  + \Gradh \vuh :  \tbS },
 \\&
  R_{\vt} =  \kappa \intTOB{ \frac{\vth^2\vthout - \tvt^3}{\tvt^2\vth\vthout } \abs{\Grad \tvt}^2  +\vth \Div\frac{\Grad \tvt}{\tvt} +  \frac{2\tvt -\vthout}{\vth \vthout} \Grad \tvt \cdot \Gradedge \vth}.
\end{align*}
Note that  $\tbS=\bS(\Grad \tvu)$ and $\bS_h =  \bS(\Gradh \vuh) = 2\mu \Dhuh +  \lambda\Divh \vuh \I$.

Next, we analyse the right hand side of \eqref{REI}. Recalling the consistency error \eqref{ecs} and thanks to the triangular inequality we estimate $R_C$ by
\begin{align*}
\abs{R_C} &\leq \Abs{e_{\vr} \left( \frac12\abs{\tvu}^2 - \frac{\pd \Hc(\tvr,\tvt)}{\pd \vr} , \TS, h,\tau\right)} + \abs{ e_{\vu} (\tvu, \TS, h,\tau)}+\abs{e_s(\tvt, \TS, h,\tau)}
\\& \aleq \TS + h^{(1-\eps)/2} +    h^{(1+\eps)/2}
\hspace{1cm} \mbox{for any } \ \eps\in(-1,1).
\end{align*}
Using Young's inequality, the uniform bounds \eqref{N} and the relation \eqref{EN} we estimate the rest terms
\begin{equation}\label{rhsRI}
\begin{aligned}
 \Abs{\sum_{i=1}^5 R_i }
 & \aleq \norm{\vuh -\tvu}_{L^2L^2}^2 + \norm{s_h -\ts}_{L^2L^2}^2
+\norm{\vrh -\tvr}_{L^2L^2}^2 + \norm{\vth -\tvt}_{L^2L^2}^2
\\& \aleq \int_0^\tau \RE{\vrh,\vuh,\vth}{\tvr,\tvu,\tvt}(t) \dt.
\end{aligned}
\end{equation}
Via H\"older's inequality, Young's inequality, the uniform bounds \eqref{N} and the relation \eqref{EN}, we can control $R_{\bS}, R_{\vt}$ as follows (see  Appendix \ref{app_res} for details):
\begin{multline}\label{REIrhs}
 \Abs{ R_{\bS}}
 \aleq   h \norm{\Gradedge \vuh}_{L^2L^2}, \
 \Abs{ R_{\vt}}
 \aleq  h + \int_0^\tau \RE{\vrh,\vuh,\vth}{\tvr,\tvu,\tvt}(t) \dt
 +\delta \norm{\Gradedge\vth-\Grad\tvt}_{L^2L^2}^2,
\end{multline}
where $h \norm{\Gradedge \vuh}_{L^2L^2} \aleq h^{(1-\eps)/2}$, cf. \eqref{N3}.
Choosing a small $\delta \in (0,{\kappa \Uvt}/ {\Ovt^2})$, we obtain the desired version of the energy inequality
\begin{equation}\label{RE04}
\begin{aligned}
&\left[ \RE{\vrh,\vuh,\vth}{\tvr,\tvu,\tvt}  \right]_0^\tau
 +\intTO{ \frac{\tvt}{ \vth} \left( 2\mu \left|\Dhuh -  \frac{\vth}{ \tvt} \bD(\tvu) \right|^2  + \lambda  \left| \Divh \vuh -  \frac{\vth}{ \tvt}\Div \tvu \right|^2 \right) }
\br
&\quad+\intTO{\frac{\kappa\tvt}{\vth \vthout }|\Gradedge\vth-\Grad\tvt|^2} \aleq \TS + h^{(1-\eps)/2} +    h^{(1+\eps)/2} + \int_0^\tau \RE{\vrh,\vuh,\vth}{\tvr,\tvu,\tvt}(t) \dt.
\end{aligned}
\end{equation}

Further, applying Gronwall's lemma, together with the error estimate caused by regular initial data
\begin{align*}
&\RE{\vrh,\vuh,\vth}{\tvr,\tvu,\tvt}(0) \aleq \intOB{ \abs{\vuh^0 -\tvu(0)}^2 + \abs{ \vth^0 -  \tvt(0)}^2 + \abs{\vrh^0 - \tvr(0)}^2 }
 \aleq h^2
\end{align*}
we obtain
\begin{align}\label{RE}
&\RE{\vrh,\vuh,\vth}{\tvr,\tvu,\tvt}(\tau)
+\intTOB{ \frac{\tvt}{ \vth} \left( 2\mu \left|\Dhuh -  \frac{\vth}{ \tvt} \bD(\tvu) \right|^2  + \lambda  \left| \Divh \vuh -  \frac{\vth}{ \tvt}\Div \tvu \right|^2 \right) }
\br
&\hspace{2cm}+ \intTO{\frac{\kappa\tvt}{\vth \vthout }|\Gradedge\vth-\Grad\tvt|^2}
\aleq \TS + h^{(1-\eps)/2} + h^{(1+\eps)/2}
\end{align}
for any $\eps \in (-1,1)$, which proves \eqref{CR1}.
Finally, \eqref{CR2} is a consequence of the uniform bounds of density and temperature, \eqref{CR1} and \eqref{EN}.
\end{proof}

\begin{Remark}
From the consistency proof and the error estimates proof above, we know that  the error of order $h^{(1-\eps)/2}$ in \eqref{RE} is indeed $h\norm{\Gradedge \vuh}_{L^2L^2}$.

Inspired by \cite{FLM18_brenner}, let us add some additional artificial viscosity in the momentum equation \eqref{scheme_M}, i.e.
\begin{align*}
\intO{ D_t  (\vrh  \vuh ) \cdot \bfphi_h } & - \intfacesint{ \Fup  (\vrh  \vuh ,\vuh ) \cdot \jump{\bfphi_h}   }+  \intO{ (\bS_h -p \I) :  \Gradh \bfphi_h }
\br
&
= - h^{\alpha} \intO{ \Gradedge \vuh  : \Gradedge (\bfphi_h) } ,
\quad \alpha > 0.
 \end{align*}
Then we obtain one more dissipation in the energy balance, i.e.
\begin{align*}
	 &D_t \intO{ \left(\frac{1}{2}  \vrh  |\vuh |^2 +  c_v \vrh  \vth \right) }
	 +  h^\eps \intfacesint{  \avs{ \vrh  }  \abs{\jump{\vuh}}^2 }
	  + \ha \intO{|\Gradedge \vuh|^2}
	\br
	&\hspace{3cm}+ \frac{\TS}{2} \intO{ \vrh^\triangleleft|D_t \vuh |^2  }
	+ \frac12 \intfacesint{ \vrh^{\rm up} |\avs{\vuh } \cdot \vc{n} | \abs{\jump{\vuh}}^2   }
	= 0,
	\end{align*}
yielding a priori bound on the velocity gradient
\begin{align*}
\norm{\Gradedge \vuh}_{L^2L^2}^2 \aleq h^{-\alpha}
\end{align*}
and a new consistency error in momentum equality, i.e. $\abs{e_{\vm}^{new}} \aleq \abs{e_{\vm}} + h^{\alpha}$.

Letting $\alpha < \eps+1$ we obtain a better estimate of $\norm{\Gradedge \vuh}_{L^2L^2}$ than \eqref{N3}.
Hence, we obtain a new formula of $A$, i.e. $A=  \min\{1-\alpha/2, \alpha, (1+\eps)/2\}$. This yields the optimal rate $A=2/3$ by choosing $\alpha = 2/3, \eps = 1/3$.
\end{Remark}

\appendix

\section{Local existence of strong solution}

The goal of this section is to prove the existence of local-in-time highly-regular strong solutions for the Navier--Stokes--Fourier system \eqref{PDE_nsf}. We pursue the same ideas developed by Kawashima \cite{Kaw} and by Breit, Feireisl, Hofmanov\'{a} \cite{BreFeiHof}: equations \eqref{PDE_nsf} can be seen as a symmetric hyperbolic-parabolic system and therefore, all the necessary bounds can be deduced from the Moser-type calculus inequalities. We can consider more general pressure $p=p(\vr, \vt)$ and internal energy $e=e(\vr, \vt)$ to be smooth functions of $\vr$ and $\vt$ satisfying 
\begin{equation} \label{pressure and internal energy}
\frac{\partial p(\vr, \vt)}{\partial \vr} \geq \Un{p}_{\vr}>0, \quad \frac{\partial e(\vr, \vt)}{\partial \vt}\geq \Un{e}_{\vt}>0.
\end{equation}

Note that condition \eqref{pressure and internal energy} implies that the thermodynamic stability condition holds. In particular, for the perfect gas we have $\frac{\partial p(\vr, \vt)}{\partial \vr}=\vt \geq \Un{\vt} >0$ and $\frac{\partial e(\vr, \vt)}{\partial \vt}=c_v = 1/(\gamma-1) > 0$.
We can now state our main result.

\begin{Theorem} \label{local ex}
	Let $k\geq 4$ and let the initial data \eqref{INI} be such that
	\begin{equation*}
	\vr_0 >0, \quad (\vr_0, \vu_0, \vt_0) \in W^{k,2}(\Of; \mathbb{R}^{1+d+1}).
	\end{equation*}
	Then there exists $T^* > 0$ and a unique strong solution $(\vr, \vu,\vt)$ to problem \eqref{PDE_nsf}, \eqref{INI} such that
	\begin{equation*}
	(\vr, \vu, \vt) \in C([0,T^*]; W^{k,2}(\Of; \mathbb{R}^{1+d+1})).
	\end{equation*}
	Moreover, the following estimate
	\begin{equation} \label{estimate ss}
	\begin{aligned}
	\sup_{t\in [0,T^*]} \Big(\| (\vr, \vu, \vt)(t) \|_{W^{k,2}} +\| \pd_{t} \vr(t) \|_{W^{k-1,2}} + \| (\pd_{t} \vu, \pd_{t} \vt)(t) \|_{W^{k-2,2}} \Big)& \\
	+\sup_{t\in [0,T^*]} \Big(\| \pd_{t}^2 \vr(t) \|_{W^{k-3,2}} + \| (\pd_{t}^2 \vu, \pd_{t}^2 \vt)(t) \|_{W^{k-4,2}} \Big)  &\leq C \| (\vr_0, \vu_0, \vt_0)\|_{W^{k,2}}
	\end{aligned}
	\end{equation}
	holds, where the positive constant $C$ depends solely on $T^*,  \ \inf \vr_0, \ \Un{e}_{\vt},$ and $\| (\vr_0, \vu_0, \vt_0)\|_{W^{k,2}}$.
\end{Theorem}

Before proving Theorem \ref{local ex}, we briefly recall some estimates which are a straightforward consequence of the Moser-type calculus, see e.g. \cite[Lemmas 2.3, 2.4 and 2.5]{Kaw}.

\begin{Lemma}\hspace{1pt}
	\begin{itemize}
		\item[(i)] Let $k \geq 2$ and $0 \leq m \leq k$. If $u \in W^{k,2}(\Of)$ and $v \in W^{m,2}(\Of)$, then $uv \in W^{m,2}(\Of)$ and
		\begin{equation} \label{moser 1}
		\| uv \|_{W^{m,2}} \leq c \| u \|_{W^{k,2}} \| v\|_{W^{m,2}}.
		\end{equation}
		\item[(ii)] Let $k \geq 1$ and let $F=F(u)$ be a $C^{\infty}$- function of $u$. If $u \in (W^{k,2} \cap L^{\infty})(\Of)$ then $\Grad F(u) \in W^{k-1,2}(\Of)$ and
		\begin{equation} \label{moser 2}
		\| \Grad F(u) \|_{W^{k-1, 2}} \leq c M (1+ \| u\|_{L^{\infty}})^{k-1} \| \Grad u\|_{W^{k-1,2}},
		\end{equation}
		with $M = \sum_{|\alpha|=1}^{k} \sup_{u} |\partial_u^{\alpha} F(u)|$ ($\sup_{u}$ is taken over all $u$ such that $|u| \leq \| u\|_{L^{\infty}}$).
		\item[(iii)] Let $k \geq 4$ and $1 \leq m \leq k$. If $u \in (W^{k,2} \cap L^{\infty})(\Of)$ and $v \in W^{m-1,2}(\Of)$, then
		\begin{equation} \label{commutator}
		[\partial_x^{\alpha}, u] v := \partial_x^{\alpha} (uv) - u \partial_x^{\alpha} v \in L^2(\Of),
		\end{equation}
		and
		\begin{equation} \label{moser 3}
		\sum_{|\alpha|=0}^{m} \| [\partial_x^{\alpha}, u] v\|_{L^2} \leq c \| \Grad u \|_{W^{k-1,2}} \| v \|_{W^{m-1, 2}}.
		\end{equation}
	\end{itemize}
\end{Lemma}

\begin{proof}[Proof of Theorem \ref{local ex}]
	As shown by Kawashima and Shizuta \cite[Section 4]{KawShi}, the Navier--Stokes--Fourier system \eqref{PDE_nsf} can be seen as a symmetric hyperbolic-parabolic system, i.e. as a coupled system of a hyperbolic equation for $\vr$ and a symmetric strongly parabolic system for $(\vu, \vt)$. More precisely, system \eqref{PDE_nsf} can be written as
	\begin{align}
	\pd_{t} \vr + \vu \cdot \Grad \vr &= f_1(\vr, \Grad \vu ), \label{continuity equation 1}\\
	\vr \pd_{t} \vu - \big[  \mu \Delta_x \vu + (\mu+\lambda) \Div \Grad^{\top} \vu \big] &= \bm{f}_2(\vr,\vu,\vt, \Grad \vr, \Grad \vu, \Grad \vt)   \label{balance of momentum 1}\\
	\vr \partial_{\vt}e \  \pd_{t} \vt - \kappa \Delta_x \vt &= f_3(\vr, \vu, \vt, \Grad \vu, \Grad \vt), \label{balance of internal energy 1}
	\end{align}
	with
	\begin{align*}
	f_1 (\vr, \Grad \vu ) :=& - \vr \Div \vu, \\
	\bm{f}_2(\vr,\vu,\vt, \Grad \vr, \Grad \vu, \Grad \vt)  :=& - \big(\partial_{\vr} p \Grad \vr+ \partial_{\vt} p \Grad \vt + \vr \vu \cdot \Grad \vu\big), \\
	f_3 (\vr, \vu, \vt, \Grad \vu, \Grad \vt) :=& \ \frac{\mu}{2}|\Grad \vu + \Grad^{\top} \vu|^2 + \lambda (\Div \vu)^2 - \big( \vt \partial_{\vt} p \  \Div \vu + \vr \partial_{\vt}e \  \vu \cdot \Grad \vt \big).
	\end{align*}
	For any $T >0$, we define the Banach space $X(T):= X_{\vr}(T) \times X_{\vu, \vt} (T)$ with
	\begin{align*}
	X_{\vr}(T):=& \ C([0,T]; W^{k,2}(\Of)) \cap C^1([0,T]; W^{k-1,2}(\Of)), \\
	X_{\vu, \vt} (T) :=&\  C ([0,T]; W^{k,2}(\Of; \mathbb{R}^{d+1})) \cap C^1([0,T]; W^{k-2,2}(\Of; \mathbb{R}^{d+1})) \cap \\
	& \ L^2(0,T; W^{k+1,2}(\Of; \mathbb{R}^{d+1})) \cap W^{1,2} (0,T; W^{k-1,2}(\Of; \mathbb{R}^{d+1})),
	\end{align*}
	and the set
	\begin{equation*}
	\mathcal{R}(T, M, M_t) := \left\{ (\vr, \vu,\vt) \in X(T)\ \Big| \ \begin{aligned}
		\sup_{t\in [0,T]} \| (\vr, \vu,\vt)(t) \|_{W^{k,2}}+  \| (\vu,\vt) \|_{L^2W^{k+1,2}}  &\leq M, \\
		\| (\pd_{t}\vr, \pd_{t} \vu, \pd_{t} \vt) (t) \|_{L^2W^{k-1,2}} &\leq M_t,
	\end{aligned} \right\}.
	\end{equation*}
	
	Let us now consider the map
	\begin{equation} \label{map G}
	(\vr, \vu, \vt) \mapsto G(\vr, \vu, \vt) = (\widetilde{\vr}, \widetilde{\vu}, \widetilde{\vt}),
	\end{equation}
	where the trio $(\widetilde{\vr}, \widetilde{\vu}, \widetilde{\vt})$ satisfies the linearized system associated to \eqref{continuity equation 1}--\eqref{balance of internal energy 1},
	\begin{align}
	\pd_{t} \widetilde{\vr} + \vu \cdot \Grad \widetilde{\vr} &= f_1(\vr, \Grad \vu ), \label{linearized 1}\\
	\vr \pd_{t} \widetilde{\vu} - \big[  \mu \Delta_x \widetilde{\vu} + (\mu+\lambda) \Div \Grad^{\top} \widetilde{\vu} \big] &= \bm{f}_2 (\vr,\vu,\vt, \Grad \vr, \Grad \vu, \Grad \vt)  \label{linearized 2} \\
	\vr \partial_{\vt}e \ \pd_{t}\widetilde{\vt} - \kappa \Delta_x \widetilde{\vt} &= f_3 (\vr, \vu, \vt, \Grad \vu, \Grad \vt). \label{linearized 3}
	\end{align}
	
	\textbf{Step 1.} Our first goal is to determine $T=T_0, M, M_t >0$ in such a way that $G: \mathcal{R}(T_0, M, M_t) \rightarrow \mathcal{R}(T_0, M, M_t)$ is well-defined. Let $(\vr, \vu, \vt ) \in \mathcal{R}(T_0, M, M_t)$ be fixed. We can deduce in particular that
	\begin{equation} \label{f}
		(f_1, \bm{f}_2, f_3) \in L^{\infty}(0,T_0; W^{k-1,2}(\Of; \mathbb{R}^{1+d+1})), \quad f_1 \in L^2(0,T_0; W^{k,2}(\Of)),
	\end{equation}
	and that there exists a constant $c_1=c_1(M)>0$ such that
	\begin{equation} \label{bound f}
		\sup_{t\in [0,T_0]}\| (f_1, \bm{f}_2, f_3)(t)\|_{W^{k-1,2}} \leq c_1 M.
	\end{equation}
	The existence of a solution $(\widetilde{\vr}, \widetilde{\vu}, \widetilde{\vt}) \in X(T_0)$ for any $T_0 >0$ of the linearized problem follows from already existing results that can be found in literature, see e.g. Kato \cite[Theorem I]{Kat}. We point out that the strict positivity of the initial density $\vr_0$ is sufficient to guarantee that
	\begin{equation*}
	\vr \geq \Un{\vr}>0.
	\end{equation*}
	In order to prove that $(\widetilde{\vr}, \widetilde{\vu}, \widetilde{\vt}) \in \mathcal{R}(T_0, M, M_t)$, we need to deduce some \textit{energy estimates}. Applying $\partial_x^{\alpha}$ with $|\alpha| \leq k$ to equations \eqref{linearized 1}--\eqref{linearized 3}, we obtain
	\begin{align}
	\pd_{t}\partial_x^{\alpha} \widetilde{\vr} + \vu \cdot \Grad \partial^{\alpha}_x \widetilde{\vr} &= F_1^{\alpha}, \label{eq 1} \\
	\vr \  \pd_{t} \partial_x^{\alpha}\widetilde{\vu} - \big[  \mu \Delta_x \partial_x^{\alpha}\widetilde{\vu} + (\mu+\lambda) \Div \Grad^{\top} \partial_x^{\alpha} \widetilde{\vu} \big] &= \bm{F}_2^{\alpha}, \label{eq 2}\\
	\vr \partial_{\vt}e \ \pd_{t} \partial_x^{\alpha}\widetilde{\vt} - \kappa \Delta_x \partial_x^{\alpha} \widetilde{\vt} &= F_3^{\alpha}, \label{eq 3}
	\end{align}
	with
	\begin{align*}
	F_1^{\alpha} &:= \partial_x^{\alpha} f_1 - \big[ \partial_x^{\alpha}, \vu \big] \Grad \widetilde{\vr}, \\
	\bm{F}_2^{\alpha} &:= \vr \ \partial_x^{\alpha} \left( \frac{\bm{f}_2}{\vr}  \right) + \vr \left(  \left[ \partial_x^{\alpha}, \frac{\mu}{\vr} \right] \Delta_x \widetilde{\vu}+ \left[ \partial_x^{\alpha}, \frac{\mu+\lambda}{\vr} \right] \Div \Grad^{\top} \widetilde{\vu} \right), \\
	F_3^{\alpha} &:= \vr \partial_{\vt}e \ \partial_x^{\alpha} \left( \frac{f_3}{\vr \partial_{\vt}e}  \right)  + \vr \partial_{\vt}e \left[ \partial_x^{\alpha}, \frac{\kappa}{\vr \partial_{\vt}e} \right] \Delta_x \widetilde{\vt};
	\end{align*}
	the commutator $\left[\partial_x^{\alpha}, \ \cdot\  \right] \ \cdot $ is defined as in \eqref{commutator}. We now multiply equations \eqref{eq 1}, \eqref{eq 2} and \eqref{eq 3} by $\partial_x^{\alpha} \widetilde{\vr}$, $\partial_x^{\alpha}\widetilde{\vu}$ and $\partial_x^{\alpha}\widetilde{\vt}$, respectively, integrate over $\Of$ and sum over $0 \leq |\alpha| \leq k$; from hypothesis \eqref{pressure and internal energy}, G\aa rding's inequality for strongly elliptic operators, H\"{o}lder and Young's inequalities, the Sobolev embedding $W^{k-1,2}(\Of) \hookrightarrow L^{\infty}(\Of)$, and from \eqref{moser 1}, \eqref{moser 2}, \eqref{moser 3}, \eqref{f},  we find constants $c_2 >0$, $c_3=c_3(\Un{\vr}, \Un{e}_{\vt}, M)>0$ and $c_4=c_4(\Un{\vr}, \Un{e}_{\vt}, M)>0$ such that
	\begin{align}
		\frac{\textup{d}}{\dt} \| \widetilde{\vr} \|_{W^{k,2}} &\leq c_2 \left( \| f_1 \|_{W^{k,2}}  +M\| \widetilde{\vr}\|_{W^{k,2}} \right), \label{e1} \\
		\frac{\textup{d}}{\dt}\| (\widetilde{\vu}, \widetilde{\vt}) \|_{W^{k,2}}^2+ c_3 \| (\widetilde{\vu}, \widetilde{\vt}) \|_{W^{k+1, 2}}^2 &\leq c_4 \left[ \left( 1+ \| ( \pd_{t} \vr, \pd_{t} \vt) \|_{W^{k-1,2}} \right) \| (\widetilde{\vu}, \widetilde{\vt})\|_{W^{k,2}}^2 + \| (\bm{f}_2, f_3)\|_{W^{k-1,2}}^2 \right]. \label{e2}
	\end{align}
	Integrating \eqref{e1}, \eqref{e2} over $[0,\tau]$, from H\"{o}lder's inequality and the Gronwall argument, there exist $c_5= c_5(\Un{\vr}, \Un{e}_{\vt})>0$ and $c_6= c_6(\Un{\vr}, \Un{e}_{\vt}, M)>0$ such that
	\begin{equation} \label{final estimate}
	\begin{aligned}
	&\| (\widetilde{\vr}, \widetilde{\vu}, \widetilde{\vt})(\tau)\|_{W^{k,2}}^2 + \int_{0}^{\tau} \| (\widetilde{\vu}, \widetilde{\vt})(t)\|_{W^{k+1,2}}^2 \ \dt \\
	&\leq c_5 e^{c_6 \left(\tau+ \tau^{\frac{1}{2}}M_t \right) } \left[ \| (\vr_0,\vu_0, \vt_0)\|_{W^{k,2}}^2 + c_6 \int_{0}^{\tau} \left(\tau \| f_1(t)\|_{W^{k,2}}^2+\| (\bm{f}_2, f_3)(t)\|_{W^{k-1,2}}^2\right) \ \dt \right].
	\end{aligned}
	\end{equation}
	From \eqref{moser 1}, \eqref{moser 2} and \eqref{bound f}, we can finally conclude that there exist $C_1=C_1(\Un{\vr}, \Un{e}_{\vt})>1$ and $C_2= C_2(\Un{\vr}, \Un{e}_{\vt}, M)>0$ such that
	\begin{align*}
	\| (\widetilde{\vr}, \widetilde{\vu}, \widetilde{\vt})(\tau)\|_{W^{k,2}}^2 &+  \| (\widetilde{\vu}, \widetilde{\vt})\|_{L^2W^{k+1,2}}^2 \leq C_1^2 e^{C_2 \left(\tau+ \tau^{\frac{1}{2}}M_t \right) } \left[ \| (\vr_0,\vu_0, \vt_0)\|_{W^{k,2}}^2 + C_2 M^2 (1+\tau)\tau \right].
	\end{align*}
	Moreover, directly from equations \eqref{continuity equation 1}--\eqref{balance of internal energy 1}, we deduce that there exists a constant $C_3=C_3(M)>0$ such that
	\begin{equation*}
	\int_{0}^{\tau} \|  ( \pd_t \widetilde{\vr}, \pd_t \widetilde{\vu}, \pd_t \widetilde{\vt})(t)\|_{W^{k-1,2}}^2 \ \dt \leq C_3^2 M^2(1+\tau).
	\end{equation*}
	To get that $(\widetilde{\vr}, \widetilde{\vu}, \widetilde{\vt}) \in \mathcal{R}(T_0, M, M_t)$, it is now enough to choose $M$, $M_t$ and $T_0$ such that
	\begin{align*}
		M &:= 2C_1 \|(\vr_0, \vu_0, \vt_0)\|_{W^{k,2}}, \\
		M_t &:= 4 C_1 C_3  \|(\vr_0, \vu_0, \vt_0)\|_{W^{k,2}}, \\
		e^{C_2 \left(T_0+ T_0^{\frac{1}{2}}M_t \right) } \leq 2, \quad C_2 M^2& (1+T_0)T_0 \leq \|(\vr_0, \vu_0, \vt_0)\|_{W^{k,2}}^2, \quad T_0 \leq 3.
	\end{align*}
	
	\textbf{Step 2.} Let us now consider the approximation sequence $\{ (\vr^n, \vu^n, \vt^n)  \}_{n\in \mathbb{N}}$ such that
	\begin{align}
		(\vr^0, \vu^0, \vt^0) &= (\vr_0, \vu_0, \vt_0), \label{first step}\\
		(\vr^n, \vu^n, \vt^n) &= G(\vr^{n-1}, \vu^{n-1}, \vt^{n-1}), \quad n\geq 1,  \label{n step}
	\end{align}
	and take the difference of two subsequent solutions:
	\begin{align*}
	\pd_{t}(\vr^{n+1}- \vr^n) + \vu^n \cdot \Grad (\vr^{n+1}- \vr^n) &= f_1^n, \\
	\vr^n \pd_{t} (\vu^{n+1}- \vu^n) - \mu \Delta_x (\vu^{n+1}- \vu^n) - (\mu+\lambda) \Div \Grad^{\top} (\vu^{n+1}- \vu^n ) &= \bm{f}_2^n, \\
	\vr^n \partial_{\vt}e^n \  \pd_{t} (\vt^{n+1}- \vt^n) - \kappa \Delta_x (\vt^{n+1}- \vt^n)&= f_3^n,
	\end{align*}
	where, from the fact that $(\vr^n, \vu^n, \vt^n) \in \mathcal{R}(T_0, M, M_t)$ for any $n \in \mathbb{N}$, it is not difficult to show that there exists a constant $C=C(\Un{\vr}, \Un{e}_{\vt}, \|(\vr_0, \vu_0, \vt_0)\|_{W^{k,2}})>0$ such that
	\begin{align*}
	\| f_1^n \|_{W^{k-1,2}} &\leq C \left( \| \vr^n -\vr^{n-1}\|_{W^{k-1, 2}} + \| (\vu^n-\vu^{n-1}, \vt^n-\vt^{n-1})\|_{W^{k,2}} \right), \\
	\| (\bm{f}_2^n, f_3^n) \|_{W^{k-2,2}} &\leq C\| (\vr^n-\vr^{n-1}, \vu^n-\vu^{n-1}, \vt^n-\vt^{n-1})\|_{W^{k-1,2}}.
	\end{align*}
	Therefore, the trio $(\vr^{n+1}- \vr^n, \vu^{n+1}- \vu^n, \vt^{n+1}- \vt^n)$ satisfies the linearized equations \eqref{linearized 1}--\eqref{linearized 3} with $(f_1, \bm{f}_2, f_3)$ replaced by $(f_1^n, \bm{f}_2^n, f_3^n)$ and  $(\vr^{n+1}- \vr^n, \vu^{n+1}- \vu^n, \vt^{n+1}- \vt^n)(0)= \textbf{0}$; hence an analogous of \eqref{final estimate} holds with $k$ replaced by $k-1$:
	\begin{align*}
	\sup_{t\in [0,\tau]} &\| (\vr^{n+1} - \vr^n, \vu^{n+1}- \vu^n, \vt^{n+1}- \vt^n)(t)\|_{W^{k-1,2}}^2 + \int_{0}^{\tau} \| (\vu^{n+1}- \vu^n, \vt^{n+1}- \vt^n)(t)\|_{W^{k,2}} \dt \\
	&\leq C_4 (1+\tau) \tau e^{C_4 \left(\tau+ M_t \tau^{\frac{1}{2}} \right)} \Big[ \sup_{t\in [0,\tau]} \| (\vr^{n} - \vr^{n-1}, \vu^{n}- \vu^{n-1}, \vt^n- \vt^{n-1})(t)\|_{W^{k-1,2}}^2 \\
	&\hspace{4.5cm}+ \int_{0}^{\tau} \| (\vu^n- \vu^{n-1}, \vt^n- \vt^{n-1})(t)\|_{W^{k,2}} \dt \Big],
	\end{align*}
	where $C_4=C_4(\Un{\vr}, \Un{e}_{\vt}, \|(\vr_0, \vu_0, \vt_0)\|_{W^{k,2}})>0$ is a constant independent of $n$. Taking $0<T_1 \leq T_0$ such that
	\begin{equation*}
	C_4T_1(1+T_1) e^{C_4(T_1+M_tT_1^{\frac{1}{2}})} <1,
	\end{equation*}
	we get that $\{ (\vr^n, \vu^n, \vt^n)  \}_{n\in \mathbb{N}}$ is a Cauchy sequence in $C([0,T_1]; W^{k-1,2}(\Of; \mathbb{R}^{1+d+1}))$. Consequently, there exists $(\vr, \vu, \vt)$ such that
	\begin{equation*}
	(\vr^n, \vu^n, \vt^n) \rightarrow (\vr, \vu, \vt)  \quad \mbox{in } C([0,T_1]; W^{k-1,2}(\Of; \mathbb{R}^{1+d+1})).
	\end{equation*}
	Furthermore, $\{ (\vr^n, \vu^n, \vt^n)  \}_{n\in \mathbb{N}} \subset \mathcal{R}(T_0, M, M_t) \subset \mathcal{R}(T_1, M, M_t)$ is uniformly bounded and therefore there exists a subsequence, not relabeled, such that
	\begin{align*}
		(\vr^n, \vu^n, \vt^n) \overset{*}{\rightharpoonup} (\vr, \vu, \vt)  &\quad \mbox{in } L^{\infty}(0,T_1; W^{k,2}(\Of; \mathbb{R}^{1+d+1})).
	\end{align*}
	Moreover, we obtain that
	\begin{equation*}
		(\pd_{t} \vr, \pd_{t} \vu, \pd_{t} \vt) \in L^{\infty}(0,T_1; W^{k-1,2}(\mathbb{T}^d) \times W^{k-2,2}(\mathbb{T}^d; \mathbb{R}^{d+1})).
	\end{equation*}
	
	\textbf{Step 3.} Finally, to show that $(\vr, \vu, \vt) \in C([0,T_1]; W^{k,2}(\mathbb{T}^d; \mathbb{R}^{1+d+1}))$, it is enough to consider a family of mollifiers $\{ \phi_{\delta} \}_{\delta>0}$ and define $(\vr_{\delta}, \vu_{\delta}, \vt_{\delta}):=(\vr * \phi_{\delta}, \vu * \phi_{\delta}, \vt * \phi_{\delta} )$. Clearly, $(\vr_{\delta}, \vu_{\delta}, \vt_{\delta}) \in C^0 ([0,T_1]; W^{k-1,2}(\mathbb{T}^d))$ for any $\delta>0$. Applying $\phi_{\delta} \ *$ to system \eqref{linearized 1}--\eqref{linearized 3} and proceeding as in the previous step, the differences $(\vr_{\delta}-\vr_{\delta'}, \vu_{\delta}- \vu_{\delta'}, \vt_{\delta}-\vt_{\delta'})$ satisfy
	the linearized equations \eqref{linearized 1}--\eqref{linearized 3} with $(f_1, \bm{f}_2, f_3)$ replaced by some $(f_1^{\delta}, \bm{f}_2^{\delta}, f_3^{\delta})$. Therefore, an analogous of \eqref{final estimate} holds and	it is easy to check that $\{ (\vr_{\delta}, \vu_{\delta}, \vt_{\delta}) \}_{\delta>0}$ is a Cauchy sequence in $C([0,T_1]; W^{k,2}(\mathbb{T}^d; \mathbb{R}^{1+d+1}))$. Consequently,
	\begin{equation*}
	(\pd_{t} \vr, \pd_{t} \vu, \pd_{t} \vt) \in C([0,T_1]; W^{k-1,2}(\mathbb{T}^d) \times W^{k-2, 2}(\mathbb{T}^d; \mathbb{R}^{d+1})),
	\end{equation*}
	and, differentiating \eqref{continuity equation 1}--\eqref{balance of internal energy 1} in time, we can deduce that
	\begin{equation*}
	(\pd_{t}^2 \vr, \pd_{t}^2 \vu, \pd_{t}^2 \vt) \in C([0,T_1]; W^{k-3,2}(\mathbb{T}^d) \times W^{k-4, 2}(\mathbb{T}^d; \mathbb{R}^{d+1})).
	\end{equation*}
	To conclude the proof, it is sufficient to choose $T^*:=T_1$.
\end{proof}

\section{Boundedness of the entropy Hessian matrix}\label{sec-EH}
In this section we derive the lower and upper bounds of the eigenvalues of the entropy Hessian matrix $\nabla_{(\vr,p)}^2(-\vr s)$ with $(\vr, \vt)\in(0,\infty)^2$.
Denote
\begin{equation*}
f(\lambda) =\Big|\nabla_{(\vr,p)}^2(-\vr s) - \lambda \mathbb{I} \Big| = \lambda^2 - \frac{c_v + (1+c_v)\vt^2}{\vr \vt^2} \lambda + \frac{c_v}{\vr^2 \vt^2}
\end{equation*}
and let $\lambda_1 \leq \lambda_2$ be the roots of $f(\lambda)=0$.

It is easy to check that
\begin{equation*}
f(0) =  \frac{c_v}{\vr^2 \vt^2} >0
\end{equation*}
and
\begin{align*}
& f\left(\frac{1}{\alpha\vr}\right) =  \frac{ [1-\alpha(1+c_v)]\vt^2 + \alpha(\alpha-1)c_v}{\alpha^2\vr^2 \vt^2} ,\quad
 f\left(\frac{1}{\beta\vr\vt^2}\right) =   \frac{[c_v\beta - (1+c_v)]\beta\vt^2 + (1-c_v\beta)}{\beta^2\vr^2 \vt^4}.
\end{align*}
Hence, we take $\alpha = 2+c_v,\,  \beta = \frac{2+c_v}{c_v}$ and obtain
\begin{align*}
1-\alpha(1+c_v) < 0, \quad & f\left(\frac{1}{\alpha\vr}\right) \geq  \frac{\Big(1- (2+c_v) \alpha + \alpha^2 \Big) c_v }{\alpha^2\vr^2 \vt^2} > 0 & \mbox{if} \ \vt^2 \leq c_v,\\
c_v\beta - (1+c_v) > 0, \quad & f\left(\frac{1}{\beta\vr\vt^2}\right) >   \frac{1+ \Big(c_v\beta - (2+c_v) \Big)c_v\beta}{\beta^2\vr^2 \vt^4} > 0 & \mbox{if} \ \vt^2 > c_v.
\end{align*}
Hence, it holds
\begin{equation*}
0 < \min\left( \frac{1}{(2+c_v)\vr}, \frac{c_v}{(2+c_v)\vr\vt^2}\right) < \lambda_1 < \frac{c_v + (1+c_v)\vt^2}{2\vr \vt^2}  < \lambda_2 < \frac{c_v + (1+c_v)\vt^2}{\vr \vt^2}.
\end{equation*}

\section{Consistency proof}\label{app_cs}
In order to prove the consistency, we integrate  \eqref{scheme_D}, \eqref{scheme_M}, and \eqref{eq_entropy_stability}  from $t=0$ to $t=t^{n+1}$ and choose $\phi_h = \Piq \phi$ as the test function for any $\phi \in L^2(0,T;W^{2,\infty}( \tor))$, $\pd_t \phi, \pd_t^2 \phi \in L^\infty((0,T)\times\tor)$. We estimate all resulting terms in six steps.
\subsection*{Step 1 -- time derivative terms}
 Let $r_h$ stand for $\vrh,$ $\vrh \vuh$ or $\vrh s_h$.
Recalling \cite[(2.17)]{FLS} and the estimates \eqref{AS}, \eqref{ESen},  and \eqref{N1} we know that
\begin{equation*}
\begin{aligned}
& \left| \left[\intTd{ r_h \phi}\right]_{t=0}^{\tau}
-  \int_0^{t^{n+1}} \intTdB{ D_t r_h(t) \Piq \phi(t) + r_h(t) \pd_t \phi (t) } \dt \right|
\\& \aleq   \TS (\norm{   \pd_{t}^2\phi  }_{L^\infty L^\infty}  + \norm{ \pd_{t}\phi   }_{L^\infty L^\infty} ) \norm{r_h}_{L^\infty L^1}
\aleq \TS
\end{aligned}
\end{equation*}
for any $\tau \in [t_n, t_{n+1}).$

In what follows we shall also need some standard interpolation inequalities used  in the numerical analysis of finite volume methods, see, e.g. \cite{EyGaHe,GallouetMAC,HS_NSF}. For completeness, we list them below.

For $\phi \in W^{2,\infty}(\tor)$ we have the following estimates,
\begin{align}
{\abs{  \jump{ \Piq  \phi  }}}
 \aleq h \| \Grad\phi \|_{L^\infty(\tor;\mathbb R^d)}, \quad
 |\Piq \phi - \avs{ \Piq  \phi} | \aleq h \| \Grad\phi \|_{L^\infty(\tor;\mathbb R^d)} \  \mbox{ for } \sigma \in \facesint \label{piqj}
\end{align}
and for all $1\leq p \leq \infty$, 
\begin{equation} \label{n4c2}
\begin{aligned}
&\norm{ \Grad \phi - \Gradedge \big(\Piq \phi\big)  }_{L^p}  \aleq h \| \Grad^2\phi \|_{L^\infty(\tor)}, \quad
\norm{ \Grad \phi - \Gradh \big(\Piq \phi\big)  }_{L^p}  \aleq h \| \Grad^2\phi \|_{L^\infty(\tor)}, \\
 &\norm{\phi - \Piq  \phi }_{L^p(\tor)} \aleq h \| \Grad\phi \|_{L^\infty(\tor;\mathbb R^d)}.
\end{aligned}
\end{equation}
Moreover, the discrete version of product rule holds
\begin{align}
\jump{f_h g_h}=\jump{f_h}\avs{g_h}+\avs{f_h}\jump{g_h}, \quad f_h, g_h \in Q_h. \label{product_rule}
\end{align}

\subsection*{Step 2 -- convective terms}
To deal with the convective terms, it is convenient to recall \cite[Lemma 2.5]{FLMS_FVNSF},
\begin{align*}
 \intn \intO{  r_h \vuh \cdot \Grad \phi  } \dt  -  \intn \intfacesint{ \Fup  [ r_h,\vuh ] \jump{  \Piq \phi}}  \dt =\sum_{j=1}^4 E_j(r_h)
\end{align*}
with
\begin{equation*}
\begin{aligned}
 E_1(r_h) &=  \frac12  \intn \intfacesint{ |\avs{\vuh} \cdot \vc{n}| \jump{r_h }  \jump{ \Piq \phi}  }  \dt,
\\ E_2(r_h) & = \frac14   \intn\intfacesint{ \jump{\vuh} \cdot \vc{n}   \jump{r_h }  \jump{  \Piq \phi}  }  \dt,
\\ E_3(r_h) & =  \intn\intO{  r_h \vuh \cdot \Big(\Grad \phi - \Gradh  \big( \Piq \phi\big) \Big) } \dt,
\\ E_4(r_h) &=
 \muh   \intn\intfacesint{   \jump{r_h }  \jump{  \Piq \phi}  }  \dt
 = {\cblue -} h^{1+\eps} \intn\intO{  r_h \Delta_h \Piq\phi } \dt.
\end{aligned}
\end{equation*}

\subsubsection*{Error terms $E_1(r_h)$}
Directly recalling the proof of  \cite[Theorem 11.2]{FeLMMiSh} we have
\begin{align*}
 \abs{E_1(r_h)}  = & \left| h \sumj \intn\intO{ r_h  \cdot \left( \Laphj{\Piq \phi} \Piq \abs{\avs{\ujh} } +  (\Piq \pdedgej{\Piq \phi})  \pdmeshj \abs{\avs{\ujh} }   \right) }   \right|
 \\
 & \aleq  h   \norm{r_h}_{L^2L^2}  \left(  \norm{\Grad^2\phi}_{L^\infty L^\infty}    \norm{\vuh}_{L^2L^2} + \norm{\Grad\phi}_{L^\infty L^\infty}  \norm{\Gradh \vuh}_{L^2L^2}  \right) \aleq h
\end{align*}
for $r_h = \vrh, \vrh \uih \ \mbox{or} \ \vrh s_h$, where $\Laphj = \pdmeshj\pdedgej$ and $\pdedgej, \pdmeshj $ are the $j^{\rm th}$-component of $\Gradedge$ and $\Gradh$, respectively.

\subsubsection*{Error terms $E_2(r_h)$}
To deal with $E_2(r_h)$  we first set $r_h = \vrh$ and obtain
\begin{equation*}
\begin{aligned}
 & \abs{E_2(\vrh)} \aleq h  \norm{\Grad\phi}_{L^\infty L^\infty}
  \left(   \intn\intfacesint{  \abs{\jump{\vuh}}^2   }  \dt  \right)^{1/2}
  \left(   \intn\intfacesint{  \jump{\vrh}^2   }  \dt  \right)^{1/2}
  \aleq h^{ (3-\eps)/2} \norm{\Gradedge \vuh}_{L^2 L^2}
\end{aligned}
\end{equation*}
due to \eqref{piqj}, \eqref{N2} and H\"older's inequality.

Then,  inserting $r_h = \vrh \uih$ into $E_2(r_h)$  and taking into account \eqref{piqj}, the assumption \eqref{AS} and the estimate \eqref{N2} with \eqref{product_rule} we get
\begin{equation*}
\begin{aligned}
 & \abs{E_2(\vrh \uih)} \aleq h  \norm{\Grad\phi}_{L^\infty L^\infty}
      \intn\intfacesint{ \big| \jump{\vuh} \cdot \vc{n} \big(  \jump{\vrh} \avs{\uih }   + \avs{\vrh} \jump{\uih }  \big)  \big|}  \dt
\\&
\aleq   h    \intn\intfacesint{  |\jump{\vrh}|\, | \jump{\vuh }\cdot \vc{n}| \,  |\avs{\vuh}|   }  \dt
+h    \intn\intfacesint{ \avs{\vrh} \abs{\jump{\vuh }}^2  }  \dt
\\&
\aleq h  \left(   \intn\intfacesint{ |\avs{\vuh} |^2 }  \dt  \right)^{1/2}
 \left(  \intn\intfacesint{  \abs{\jump{\vuh }}^2 }  \dt  \right)^{1/2}
 +h^2 \norm{\Gradedge \vuh}_{L^2 L^2}^2
 \\& \aleq  h \norm{\Gradedge \vuh}_{L^2 L^2}  + h^2 \norm{\Gradedge \vuh}_{L^2 L^2}^2.
 \end{aligned}
\end{equation*}
From the Taylor expansion
\begin{equation*}
 |\jump{\vrh s_h} | = \abs{\pd_{\vr}(-\vr s)(\vv^*) \jump{\vrh} +  \pd_{p}(-\vr s)(\vv^*) \jump{p_h} } \aleq \abs{\jump{\vrh}} + \abs{\jump{p_h}}
\end{equation*}
with $\vv^*\in \co{(\vrh^{\rm in}, p_h^{\rm in})}{(\vrh^{\rm out}, p_h^{\rm out})}$. As $\jump{p_h}$ has the same bound as $\jump{\vrh}$, see  \eqref{N2}, we know that $E_2(\vrh s_h)$ has the same estimate as $E_2(\vrh)$, meaning
\begin{equation*}
 \abs{E_2(\vrh s_h)}
 \aleq h^{  (3-\eps)/2} \norm{\Gradedge \vuh}_{L^2 L^2}.
\end{equation*}

\subsubsection*{Error terms $E_3(r_h)$}
 The estimates of the third error terms are straightforward due to \eqref{n4c2}, and the uniform bounds \eqref{AS} -- \eqref{N1}. Indeed,
\begin{equation*}
\begin{aligned}
 \abs{E_3(\vrh)} & \aleq  h  \norm{\Grad^2\phi}_{L^\infty L^\infty} \intn\intO{ | \vrh \vuh |} \dt \aleq h \norm{\vuh}_{L^\infty L^2} \aleq h ,
\\ \abs{E_3(\vrh \uih)} & \aleq h  \norm{\Grad^2\phi}_{L^\infty L^\infty}  \intn\intO{ | \vrh \uih \vuh |} \dt \aleq h \norm{\vuh}_{L^\infty L^2}^2  \aleq h ,
\\ \abs{E_3(\vrh s_h)} & \aleq  h  \norm{\Grad^2\phi}_{L^\infty L^\infty}  \intn\intO{ | \vrh s_h \vuh |} \dt \aleq h \norm{\vuh}_{L^\infty L^2}  \aleq h.
\end{aligned}
\end{equation*}

\subsubsection*{Error terms $E_4(r_h)$}
Finally, we treat  the fourth error terms.
For $r_h=\vrh s_h$ the term is not present, i.e. $E_4(\vrh s_h)=0$. For $r_h=\vrh, \vrh \uih$ we have
\begin{equation*}
 \abs{E_4(r_h)} \aleq h^{1+\eps}  \norm{\Grad^2\phi}_{L^\infty L^\infty} \norm{r_h}_{L^1 L^1} \aleq h^{1+\eps}.
 \end{equation*}

Similarly as in the first step, recalling \eqref{AS} and  \eqref{N1} we realize that
\begin{align}\label{aux}
 \left| \int_\tau^{t^{n+1}}\intO{  r_h \vuh \cdot \Grad \phi  } \dt \right| \aleq  \TS\|\Grad\phi\|_{L^\infty L^\infty}\|\vuh\|_{L^\infty L^2}  \|r_h\|_{L^\infty L^2}\aleq \TS
\end{align}
for $r_h \in \{\vrh, \vrh u_{i,h}, \vrh s_h\}.$
\medskip
\noindent
 Collecting the above estimates of $E_j(r_h), j=1,\ldots, 4,$ and \eqref{aux} we get
\begin{equation}\label{CST2}
\begin{aligned}
&\Abs{ \intTO{  r_h \vuh \cdot \Grad \phi  } \dt  -  \intn\intfacesint{ \Fup  [ r_h,\vuh ] \jump{  \Piq \phi}}  \dt }
\aleq \TS + h + h^{1+\eps} + c(r_h),
\\&
c(r_h) \aleq \left\{ \begin{array}{ll}
h^{(3-\eps)/2}\norm{\Gradedge \vuh}_{L^2L^2} \aleq h^{1-\eps}  & \mbox{for } r_h=\vrh, \vrh s_h, \\
h\norm{\Gradedge \vuh}_{L^2L^2} +h^2\norm{\Gradedge \vuh}_{L^2L^2}^2 \aleq h^{(1-\eps)/2}+ h^{1-\eps} & \mbox{for } r_h=\vrh u_{i,h}.
\end{array} \right.
\end{aligned}
\end{equation}
\subsection*{Step 3 -- viscosity terms in the momentum equation}
 The interpolation error estimate (\ref{n4c2}) and the {\sl a priori} bound \eqref{N1} are enough to control  the consistency error of the viscosity terms in the momentum equation. Indeed,   we have
 \begin{equation*}
\begin{split}
&  \Abs{ \intTO{ \bS_h : \Grad  \bfphi  }  -  \intn\intO{ \bS_h : \Gradh(\Piq \bfphi) }\dt }
\\& \leq
 \Abs{ \intn\intO{ \bS_h : \big(  \Grad  \bfphi  - \Gradh (\Piq \bfphi)  \big) }\dt}
 +\Abs{ -\int_\tau^{t^{n+1}}\intO{ \bS_h : \Grad  \bfphi}}
\\&
\lesssim
h \norm{\bS_h }_{L^2 L^2}   \norm{\Grad^2\bfphi}_{L^\infty L^\infty}
 +\Abs{ \int_\tau^{t^{n+1}}\intO{ \Gradh \vuh : \Piq \bS(\Grad \bfphi)}}
\\&
\lesssim
 h  +\Abs{ \int_\tau^{t^{n+1}}\intO{ \vuh \cdot \Divh \Piq \bS(\Grad \bfphi)}}
\\&
\aleq h + (t^{n+1}-\tau)  \norm{\vuh}_{L^\infty L^2} \norm{ \bfphi}_{L^\infty W^{2,\infty}}
 \aleq h+\TS
\end{split}
\end{equation*}
where we have used the fact that $\bS_h$ is piecewise constant and the identity
\begin{equation}\label{sym}
 \bS_h :\Grad \bfphi =  \Gradh \vuh :\bS(\Grad \bfphi)
\end{equation}
due to the symmetry of $\bS_h$ and $\bS(\Grad\bfphi)$.
\subsection*{Step 4 -- pressure term in the momentum equation}
The consistency error of the pressure term in the momentum equation is controlled, thanks to the interpolation estimate \eqref{n4c2} and the uniform bound on pressure  \eqref{ESen}, and the triangular inequality.
\begin{align*}
& \bigg| \int_0^{\tau}\intO{ p_h  \Div\bfphi }\dt   -  \intn\intO{ p_h  \Divh ( \Piq  \bfphi)  }\dt \bigg|
\\& \leq \Abs{\intn \intO{ p_h  \big( \Div\bfphi-  \Divh ( \Piq  \bfphi) \big)  }} +\Abs{ \int_\tau^{t^{n+1}}\intO{ p_h  \Div\bfphi   } \dt }
\\& \leq \norm{p_h}_{L^\infty L^1} h  \norm{\Grad^2\bfphi}_{L^\infty L^\infty}
+ (t^{n+1}-\tau) \norm{p_h}_{L^\infty L^1}  \norm{\Div \bfphi}_{L^\infty L^\infty}  \aleq \TS + h.
\end{align*}

\subsection*{Step 5 -- $\kappa$-term in the entropy inequality}
 Using the product rule \eqref{product_rule}, i.e. $\jump{\frac{\Piq \phi}{\vth}}  =  \jump{\Piq \phi} \avs{\frac{1}{\vth}} + \avs{\Piq \phi} \jump{ \frac{1}{\vth} } = \jump{\Piq \phi} \avs{\frac{1}{\vth}} - \avs{\Piq \phi} \frac{\jump{\vth }}{\vthout \vth} $, we can write the consistency error of the $\kappa$-term in the following way
\begin{align*}
 &	-\int_0^{t^{n+1}} \intfacesint{\frac{\kappa}{ h } \jump{\vth}   \jump{ \frac{\Piq \phi }{\vth }} } \dt
 - \left(   \intTO{ \frac{ \kappa \phi}{ \vthout \vth} \abs{\Gradedge \vth }^2 }
-   \intTO{ \frac{\kappa}{\vth} \Gradedge \vth \cdot \Grad \phi }   \right)
\\ & =
\kappa \int_\tau^{t^{n+1}} \intO{\frac{\avs{\Piq \phi } }{\vth \vthout} \abs{\Gradedge\vth}^2   } \dt
  + \kappa \intTO{\frac{\avs{\Piq \phi}   - \phi}{\vth \vthout} \abs{\Gradedge\vth}^2   } \dt
\\& 	\quad
+ \intTO{ \frac{\kappa}{\vth} \Gradedge \vth \cdot (\Grad \phi   - \Gradedge \Piq \phi)}
+\kappa \intTO{ \left(\frac{1}{\vth} - \avs{ \frac{1}{\vth}} \right)\Gradedge \vth \cdot  \Gradedge \Piq \phi}
 \\&\quad
 -\kappa \int_\tau^{t^{n+1}} \intO{ \avs{ \frac{1}{\vth}}   \Gradedge\vth  \cdot   \Gradedge \Piq \phi   } \dt
\\ &=:   \kappa \int_\tau^{t^{n+1}} \intO{\frac{\avs{\Piq \phi } }{\vth \vthout} \abs{\Gradedge\vth}^2   } \dt
+ \sum_{i=1}^4 I_i  .
\end{align*}
Note that $ \kappa \int_\tau^{t^{n+1}} \intO{\frac{\avs{\Piq \phi } }{\vth \vthout} \abs{\Gradedge\vth}^2   } \dt  \geq 0$ contributes to the $``\geq "$ sign of \eqref{cP3}.
 The $I_i$, $i=1,\ldots,4$, terms contribute to the consistency error, and can be controlled by using H\"older's inequality, the uniform bounds \eqref{AS} and \eqref{N}, and the integration by part formula \eqref{dis_op3}, i.e.
\begin{align*}
\abs{I_1} &
\aleq  h \norm{\phi}_{L^\infty W^{1,\infty}} \norm{\Gradedge \vth}_{L^2L^2} \aleq h, \quad
 \abs{I_2}
\aleq h \norm{\phi}_{L^2 W^{2,\infty}} \norm{\Gradedge \vth}_{L^2L^2} \aleq h,
\\ \abs{I_3} &
=  \frac{h \kappa}2 \Abs{\intTO{\frac{\Gradedge \vth }{\vthout \vth}  \Gradedge \vth \cdot  \Gradedge \Piq \phi}  }
\aleq h\norm{\phi}_{L^\infty W^{1,\infty}}   \norm{\Gradedge \vth}_{L^2L^2}^2
\aleq h,
\\   \abs{I_4} &
\aleq  \Abs{ \int_\tau^{t^{n+1}} \intO{  \Gradedge\vth  \cdot   \Gradedge \Piq \phi   } \dt}
 = \Abs{ \int_\tau^{t^{n+1}} \intO{  \vth  \cdot   \Laph \Piq \phi   } \dt} \aleq \TS \norm{\phi}_{L^\infty W^{2,\infty}}\aleq \TS.
\end{align*}

\subsection*{Step 6 -- entropy production terms in the entropy inequality}
Thanks to $\vth > 0$ and $\phi \geq 0$ the entropy production terms are non-negative, i.e.
\begin{align*}
& \intn\intO{D_1 (\Piq\phi)  + D_2(\Piq\phi)+ D_3(\Piq\phi) } \dt \geq 0,
\\
&  \intn\intO{\bS_h:\Gradh \vuh \frac{\Piq\phi}{\vth} } \dt
-\intTO{ \bS_h:\Gradh \vuh \frac{ \phi}{\vth} } \dt
\\
& = \int_{\tau}^{t_{n+1}}\intO{ \bS_h:\Gradh \vuh  \frac{ \phi}{\vth} } \dt
= \int_{\tau}^{t_{n+1}}\intO{ \left(2\mu |\Dhuh|^2  + \lambda |\Divh \vuh|^2  \right) \frac{ \phi}{\vth} } \dt \geq 0,
\end{align*}
which again contribute to the $``\geq "$ sign of \eqref{cP3}.

In addition, by \eqref{piqj}, \eqref{N2} and H\"older's inequality, we estimate the residual term in the entropy balance \eqref{eq_entropy_stability}
\begin{equation*}
 \Abs{ \intn \intO{R_{s}(\Piq\phi) }}
\aleq h^{1+\eps}  \norm{\Grad \phi}_{L^\infty L^\infty} \int_0^{\tau}\intfacesint{    \left( \jump{\vrh} + \jump{p_h}\right)      }\dt
\aleq  h^{(\eps+1)/2} +   h^{\eps+1}.
\end{equation*}
Finally, collecting the estimates of the six steps above, we obtain the consistency formulation \eqref{cP1}, \eqref{cP2}, and \eqref{cP3} from \eqref{scheme_D}, \eqref{scheme_M}, and \eqref{eq_entropy_stability}, respectively. This proves Lemma \ref{lem_C1}.

\section{Relative energy inequality}\label{app_REI}
In this section we derive the relative energy inequality \eqref{REI}. More precisely, we want to prove
\begin{align}\label{AREI}
&\left[ \RE{\vrh,\vuh,\vth}{\tvr,\tvu,\tvt}  \right]_0^\tau
 +\int_0^\tau \intOB{ \frac{\tvt}{ \vth} \left( 2\mu \left|\Dhuh -  \frac{\vth}{ \tvt} \bD(\tvu) \right|^2  + \lambda  \left| \Divh \vuh -  \frac{\vth}{ \tvt}\Div \tvu \right|^2 \right) }
\br
&\hspace{2cm} +\intTO{\frac{\kappa\tvt}{\vth \vthout }|\Gradedge\vth-\Grad\tvt|^2} \leq
R_C   +\sum_{i=1}^5 R_i - R_{\bS} - R_\vt,
\end{align}
where
\begin{equation}\label{ARIS}
\begin{aligned}
&R_C = e_{\vr} \left( \frac12\abs{\tvu}^2 - \frac{\pd \Hc(\tvr,\tvt)}{\pd \vr} , \TS, h,\tau \right) +e_{\vu} (\tvu, \TS, h,\tau) - e_s(\tvt, \TS, h,\tau),
\\& R_1=
 - \intTO{  \vrh (s_h - \ts)  (\vuh-\tvu) \cdot \Grad \tvt  }  , \quad
R_2= - \int_0^\tau \intO{  \vrh (\vuh - \tvu) \otimes  (\vuh - \tvu): \Grad \tvu },
\\& R_3= \int_0^\tau \intO{   \frac{ \vrh-\tvr}{\tvr} ( \tvu - \vuh )\cdot  \Div \tbS  }
,  \quad
R_4= \intTO{ ( \tp-p_h   -  \pd_\vr \tp (\tvr - \vrh)   -  \pd_\vt \tp (\tvt - \vth)   )\Div \tvu }  ,
\\& R_5=
-\intTO{  \bigg(  (\vrh-\tvr) (s_h- \ts)  + \tvr \big(s_h- \ts - \pd_\vr \ts (\vrh -\tvr)  - \pd_\vt \ts (\vth -\tvt) \big)  \bigg)(\pd_t\tvt + \tvu \cdot \Grad \tvt ) }   ,
\\&
 R_{\bS} =
\intTOB{  \vuh \cdot  \Div \tbS  + \Gradh \vuh :  \tbS },
 \\&
 R_{\vt} =  \kappa \intTOB{ \frac{\vth^2\vthout - \tvt^3}{\tvt^2\vth\vthout } \abs{\Grad \tvt}^2  +\vth \Div\frac{\Grad \tvt}{\tvt} +  \frac{2\tvt -\vthout}{\vth \vthout} \Grad \tvt \cdot \Gradedge \vth}.
\end{aligned}
\end{equation}

\paragraph{Step 1.}  Our first goal is to obtain a relative energy inequality directly from the consistency formulation \eqref{cP} and the energy balance \eqref{energy_stability}. We reformulate the relative energy as
\begin{equation}\label{RE0}
\RE{\vrh,\vuh,\vth}{\tvr,\tvu,\tvt} = \sum_{i=1}^4 T_i,
\end{equation}
where
\begin{align*}
&T_1 = \intOB{\frac12 \vrh \abs{\vuh}^2 + \Hc(\vrh,\vuh)},
\quad
T_2 = \intO{\vrh \left( \frac12  \abs{\tvu}^2 -  \frac{\pd \Hc(\tvr,\tvt)}{\pd \vr} \right) } ,
\\ &
T_3= - \intO{ \vrh \vuh \cdot \tvu},
\quad
T_4=   \intOB{ \tvr \frac{\pd \Hc(\tvr,\tvt)}{\pd \vr} -  \Hc(\tvr,\tvt)}= \intO{ \tp} .
\end{align*}
Hence, for the first term $T_1$, by integrating \eqref{energy_stability} from $t=0$ to $t=\tau$ and subtracting \eqref{cP3} with $\phi=\tvt$, we obtain
\begin{equation}\label{RE1}
\begin{aligned}
\left[ T_1 \right]_{t=0}^\tau 
&\leq
-   \intTO{ \frac{\kappa \tvt }{ \vthout  \vth} \abs{\Gradedge \vth}^2  }
 -  \intTO{\bS_h:\Gradh \vuh \frac{ \tvt}{\vth} }
\\&
- \intTO{  \vrh s_h (\pd_t\tvt + \vuh \cdot \Grad \tvt ) }
+   \intTO{  \frac{\kappa}{\vth} \Gradedge \vth \cdot \Grad \tvt  }
- e_{s}\left(\tvt, \TS, h, \tau\right) .
\end{aligned}
\end{equation}
For the second term $T_2$, by taking $\phi = \frac12\abs{\tvu}^2 - \frac{\pd \Hc(\tvr,\tvt)}{\pd \vr} $ in \eqref{cP1} we obtain
\begin{equation} \label{RE2}
\begin{aligned}
\left[ T_2 \right]_{t=0}^\tau
& =
  \int_0^\tau \intOB{   \vrh \partial_t \tvu \cdot \tvu +  \vrh  \tvu \otimes \vuh : \Grad \tvu  } \dt   + \int_0^\tau \intO{   \vrh \ts( \partial_t \tvt +   \vuh \cdot \Grad \tvt )} \dt
\\&
 -\int_0^\tau \intO{  \frac{ \vrh}{\tvr} (\partial_t \tp +   \vuh \cdot \Grad \tp) } \dt
  +  e_\vr\left(\frac12\abs{\tvu}^2 - \frac{\pd \Hc(\tvr,\tvt)}{\pd \vr}, \TS, h,\tau\right).
 \end{aligned}
\end{equation}
For the third term $T_3$, by setting $\bfphi=\tvu$ in \eqref{cP2} we obtain
\begin{align} \label{RE4}
\left[ T_3 \right]_{t=0}^\tau  &=
- \intTOB{  \vrh \vuh \cdot \partial_t \tvu + \vrh \vuh \otimes \vuh  : \Grad \tvu }
\br
&+\intTO{(\bS_h - p_h\I):\Grad \tvu}
- e_{\vm}\left(\tvu, \TS, h,\tau\right).
\end{align}
For the fourth the term $T_4$ we simply write
\begin{equation}\label{RE5}
 \left[ T_4 \right]_{t=0}^\tau =   \intTO{ \pd_t \tp }  .
\end{equation}
Now, summing up \eqref{RE1}--\eqref{RE5} we obtain from $\left[ \RE{\vrh,\vuh,\vth}{\tvr,\tvu,\tvt}  \right]_0^\tau  = \sum_{i=1}^4 [T_i]_{t=0}^\tau$ that
\begin{align}\label{RE02}
&\left[ \RE{\vrh,\vuh,\vth}{\tvr,\tvu,\tvt}  \right]_0^\tau
  + \Xi_{l}
 \leq
R_C   
- \intTO{  \vrh (s_h-\ts)  (\vuh-\tvu) \cdot \Grad \tvt  } + \Xi_{r}
\end{align}
with
\begin{align*}
\Xi_{l} & = \intTOB{ \bS_h: \Gradh \vuh \frac{ \tvt}{\vth}  - \bS_h: \Grad \tvu }
 +   \intTOB{ \frac{\kappa \tvt }{ \vthout  \vth} \abs{\Gradedge \vth}^2  -  \frac{\kappa}{\vth} \Gradedge \vth \cdot \Grad \tvt  },
\br
\Xi_{r} & =
\intTOB{   \vrh (\tvu -\vuh) \cdot \partial_t \tvu  +   \vrh ( \tvu-\vuh) \otimes \vuh : \Grad \tvu }  + \intTO{  \vrh (\ts - s_h) (\pd_t\tvt + \tvu \cdot \Grad \tvt ) }
\br
&\quad
+ \intTOB{ \pd_t \tp  -  p_h \Div \tvu } \dt -\int_0^\tau \intO{  \frac{ \vrh}{\tvr} (\partial_t \tp +   \vuh \cdot \Grad \tp) }  .
\end{align*}
We point out that the derivative-terms of $\Xi_{r}$ are the first order derivative-terms of the exact solution $(\tvr, \tvu, \tvt)$ so that all of them could be bounded by some constants.

\paragraph{Step 2.}  The second goal is to reformulate $\Xi_{r}$ into the difference between $(\vrh, \vuh, \vth)$ and $(\tvr, \tvu, \tvt)$ so that it can be controlled by the relative energy with \eqref{EN}.

The key is to utilize that $(\tvr, \tvu, \tvt)$ is the strong solution satisfying
\begin{align*}
& \pd_t \tvr + \tvu \cdot \Grad \tvr = - \tvr \Div \tvu, \quad
 \tvr (\pd_t \tvu + \tvu \cdot \Grad \tvu)  + \Grad \tp   = \Div \tbS,
 \br
 &\partial_t \tvt + \tvu \cdot \Grad \tvt + \frac{\tvt}{c_v} \Div \tvu = \frac{1}{c_v \tvr} \left( \tbS : \Grad \tvu + \frac{\kappa \abs{\Grad \tvt}^2}{\tvt} + \vt \Div\left( \frac{\kappa \Grad \tvt}{\tvt} \right)\right) .
\end{align*}
Firstly, with
\begin{align*}
&  \intOB{    \vrh( \tvu - \vuh )\cdot \partial_t \tvu  +  \vrh ( \tvu-\vuh) \otimes \vuh : \Grad \tvu}
\\ = \ &
 \intO{    \vrh( \tvu - \vuh )\cdot \left( \partial_t \tvu  +   \tvu  \cdot \Grad \tvu \right)  }
 -    \intO{  \vrh (\vuh - \tvu) \otimes  (\vuh - \tvu): \Grad \tvu }
\\ = \ &
 \intO{   \frac{ \vrh}{\tvr}( \tvu - \vuh )\cdot \left( \Div \bS(\Grad \tvu) - \Grad \tp \right)  }
 -    \intO{  \vrh (\vuh - \tvu) \otimes  (\vuh - \tvu): \Grad \tvu }
\end{align*}
and $\intOB{  \tvu \cdot \Grad \tp +   \tp  \Div \tvu }  =0$, we obtain
\begin{align*}
\Xi_{r} & =
-  \intTO{  \vrh (\vuh - \tvu) \otimes  (\vuh - \tvu): \Grad \tvu } + \intTO{   \frac{ \vrh}{\tvr}( \tvu - \vuh )\cdot  \Div \tbS  } + \intTO{\Xi_{r,1}},
\\
\Xi_{r,1} & =   \vrh ( \ts - s_h) (\pd_t\tvt + \tvu \cdot \Grad \tvt ) + \frac{ \tvr - \vrh}{\tvr} (\pd_t \tp  + \tvu \cdot \Grad \tp ) + ( \tp-p_h)\Div \tvu .
\end{align*}
Secondly, applying the product rule for $\tp$ and $\ts$ as functions of $(\tvr,\tvt)$
\begin{align*}
& (\pd_t \tp  + \tvu \cdot \Grad \tp )
 =
  (\pd_t \tvr  + \tvu \cdot \Grad \tvr) \pd_\vr \tp+
 (\pd_t\tvt + \tvu \cdot \Grad \tvt )\pd_\vt \tp =
 -\tvr \Div \tvu  \pd_\vr \tp+
\tvr  (\pd_t\tvt + \tvu \cdot \Grad \tvt )
\br
& \vrh (s_h- \ts)
 = (\vrh-\tvr) (s_h- \ts)
+ \tvr (s_h- \ts - \pd_\vr \ts (\vrh -\tvr)  - \pd_\vt \ts (\vth -\tvt) )
+ \tvr ( \pd_\vr \ts (\vrh -\tvr)  + \pd_\vt \ts (\vth -\tvt) )
\end{align*}
we have
\begin{align*}
\Xi_{r,1} = \  &  \bigg(\tp-p_h -  \pd_\vr \tp( \tvr - \vrh) \bigg) \Div \tvu  + \bigg( \vrh ( \ts - s_h) + ( \tvr - \vrh) \bigg) (\pd_t\tvt + \tvu \cdot \Grad \tvt )
\br
= \ &  \bigg(\tp-p_h -  \pd_\vr \tp( \tvr - \vrh) -  \pd_\vt \tp (\tvt - \vth) \bigg) \Div \tvu
\br
& - \bigg( (\vrh-\tvr) (s_h- \ts)
+ \tvr (s_h- \ts - \pd_\vr \ts (\vrh -\tvr)  - \pd_\vt \ts (\vth -\tvt) ) \bigg) (\pd_t\tvt + \tvu \cdot \Grad \tvt )
\br
& -(\vth -\tvt) \bigg(\tvr \pd_\vt \ts  (\pd_t\tvt + \tvu \cdot \Grad \tvt ) + \pd_\vt \tp \Div \tvu  \bigg).
\end{align*}
With $\pd_\vt \ts = \frac{c_v}{\tvt}, \pd_\vt \tp = \tvr$, the last term above can be rewritten as
\begin{align*}
-(\vth -\tvt) \bigg(\tvr \pd_\vt \ts  (\pd_t\tvt + \tvu \cdot \Grad \tvt ) + \pd_\vt \tp \Div \tvu  \bigg) = -\frac{\vth -\tvt}{\tvt} \left( \tbS : \Grad \tvu + \frac{\kappa \abs{\Grad \tvt}^2}{\tvt} + \tvt \Div\left( \frac{\kappa \Grad \tvt}{\tvt} \right)\right).
\end{align*}

Thanks to
\begin{align*}
 \intOB{   \tvu   \cdot \Div \tbS  + \tbS: \Grad \tvu  } = 0, \quad  \intTO{\frac{\kappa \abs{\Grad \tvt}^2}{\tvt} + \tvt \Div\left( \frac{\kappa \Grad \tvt}{\tvt} \right)} = 0,
\end{align*}
finally we reformulate $\Xi_{r}$ as
\begin{align*}
\Xi_{r} & =
-  \intTO{  \vrh (\vuh - \tvu) \otimes  (\vuh - \tvu): \Grad \tvu } + \intTO{   \frac{ \vrh-\tvr}{\tvr}( \tvu - \vuh )\cdot  \Div \tbS  }
\br
& + \intTO{\bigg(\tp-p_h -  \pd_\vr \tp( \tvr - \vrh) -  \pd_\vt \tp (\tvt - \vth) \bigg) \Div \tvu}
\br
& - \intTO{\bigg( (\vrh-\tvr) (s_h- \ts)
+ \tvr (s_h- \ts - \pd_\vr \ts (\vrh -\tvr)  - \pd_\vt \ts (\vth -\tvt) ) \bigg) (\pd_t\tvt + \tvu \cdot \Grad \tvt )} - \Xi_{r,1}
\br
 \Xi_{r,1} & = \intTOB{\frac{\vth }{\tvt}  \tbS : \Grad \tvu + \vuh \cdot  \Div \tbS } +\intTO{\frac{\vth }{\tvt} \left(\frac{\kappa \abs{\Grad \tvt}^2}{\tvt} + \tvt \Div\left( \frac{\kappa \Grad \tvt}{\tvt} \right)\right)}.
\end{align*}

\paragraph{Step 3.}  Back to the relative energy inequality \eqref{RE02}, the rest is to reformulate $\Xi_{l} +  \Xi_{r,1}$.
We split $\Xi_{l} +  \Xi_{r,1}$ into the velocity-gradient-terms and temperature-gradient-terms:
\begin{equation*}
\begin{aligned}
& \Xi_{l} +  \Xi_{r,1} = T_{\bS} + T_\vt ,
\\
&T_{\bS}  =
  \intTOB{ \bS_h:\Gradh \vuh \frac{ \tvt}{\vth}  - \bS_h: \Grad \tvu + \frac{\vth }{\tvt}  \tbS : \Grad \tvu + \vuh \cdot  \Div \tbS},
\\
&T_\vt  =
  \intTO{ \frac{\kappa \tvt }{ \vthout  \vth} \abs{\Gradedge \vth}^2  -  \frac{\kappa}{\vth} \Gradedge \vth \cdot \Grad \tvt + \frac{\vth }{\tvt} \left(\frac{\kappa \abs{\Grad \tvt}^2}{\tvt} + \tvt \Div\left( \frac{\kappa \Grad \tvt}{\tvt} \right)\right)}.
\end{aligned}
\end{equation*}

\noindent{\bf Term $T_{\bS}$:} Denoting $R_{\bS} = \intTOB{  \vuh \cdot  \Div \tbS  + \bS_h :  \Grad \tvu }$, we have
\begin{equation*}
\begin{aligned}
  T_{\bS}
 = &
  \intTOB{ \frac{\tvt}{ \vth}
 \bS_h :\Gradh \vuh -  2\bS_h : \Grad \tvu + \frac{\vth}{ \tvt}   \tbS :\Grad \tvu }
 +   R_{\bS}
 \\=&
  \intTO{ \frac{\tvt}{ \vth} \left( 2\mu \left|\Dhuh -  \frac{\vth}{ \tvt} \bD(\tvu) \right|^2  + \lambda  \left| \Divh \vuh -  \frac{\vth}{ \tvt}\Div \tvu \right|^2 \right) }
 + R_{\bS}.
 \end{aligned}
 \end{equation*}

\noindent{\bf Term $T_\vt$:} We straightforward reformulate $T_\vt$ as
\begin{align*}
R_{\vt} = T_\vt -
\intTO{\frac{\kappa\tvt}{\vth \vthout }|\Gradedge\vth-\Grad\tvt|^2}.
\end{align*}
Finally, collecting the above and going back to \eqref{RE02} we obtain the relative energy inequality \eqref{AREI}.

\section{Estimates on $R_{\bS}$ and $R_{\vt}$}\label{app_res}
In this part we prove \eqref{REIrhs}, namely the following estimate.
\begin{multline*}
 \Abs{ R_{\bS}}
 \aleq   h \norm{\Gradedge \vuh}_{L^2L^2},
 \Abs{ R_{\vt}}
 \aleq  h + \int_0^\tau \RE{\vrh,\vuh,\vth}{\tvr,\tvu,\tvt}(t) \dt
 +\delta \norm{\Gradedge\vth-\Grad\tvt}_{L^2L^2}^2.
 \end{multline*}

\paragraph{Estimate of $R_{\bS}$}
Noticing the fact that $\Gradh \vuh$ is piecewise constant, we have
\begin{align*}
&R_{\bS} =  \intTOB{  \vuh \cdot  \Div \tbS  + \Gradh \vuh :  \Piq \tbS }
=  \intTO{  \vuh \cdot  (\Div \tbS  - \Divh \Piq \tbS )}
\\& = \int_0^\tau \sum_K \vuh|_K\cdot   \sum_{\sigma \in \facesK}\intSh{\left(\tbS  - \avs{\Piq \tbS}\right)\cdot \vn}\dt
\\& =- \int_0^\tau \intfacesint{ \jump{\vuh} \cdot \left(\tbS  - \avs{\Piq \tbS}\right)\cdot \vn }\dt
 \aleq h \norm{\tvu}_{L^\infty W^{2,\infty}} \int_0^\tau \intfacesint{ \abs{ \jump{\vuh} } }\dt
\\& \aleq h \left(   \int_0^\tau\intfacesint{  \abs{\jump{\vuh}}^2   }  \dt  \right)^{1/2}
  \left(   \int_0^\tau\intfacesint{  1  }  \dt  \right)^{1/2}
 = h \norm{\Gradedge \vuh}_{L^2L^2} .
\end{align*}

\paragraph{Estimate of $R_{\vt}$}
We start with the reformulation $ R_{\vt} = \sum_{i=1}^5 R_{\vt,i}$, where
\begin{equation*}
\begin{aligned}
 R_{\vt,1} = &
\kappa  \intTOB{ \vth \Div \left(\frac{ \Grad \tvt}{\tvt} \right) +  \Gradedge\vth \cdot \frac{\Grad \tvt  }{\tvt}  }  ,
\\
R_{\vt,2} = &  \kappa \intTO{ \left(  \frac{\vth}{\tvt}-  \frac{\tvt}{\vth} \right) \left(\frac{1}{\tvt}-\frac{1}{ \vthout }\right)|\Grad \tvt|^2}\dt,
\\
R_{\vt,3} = &  \kappa \intTO{\frac{1}{\vthout}\left(\frac{\vth}{\tvt}-\frac{\tvt}{\vth}\right)\left(\Grad \tvt - \Gradedge\vth\right) \cdot \Grad \tvt},
\\
R_{\vt,4} = &   \kappa \intTO{\frac{\tvt}{\vth}\left(\frac{1}{\tvt}-\frac{1}{ \vthout }\right)(\Grad \tvt-\Gradedge\vth)\cdot\Grad \tvt } ,
\\
R_{\vt,5} = &  \kappa\intTO{ \frac{1}{\tvt} \left( \frac{\vth}{\vthout} - 1\right) \Gradedge\vth\cdot \Grad \tvt  }
 =  - h\kappa  \intTO{\frac{1}{\tvt \vthout } |\Gradedge\vth|^2\Grad\tvt }.
\end{aligned}
\end{equation*}
Hence, with \eqref{dis_op3}, the assumption \eqref{AS}, the projection error \eqref{proj} and the a priori bound on temperature \eqref{N1} we have
\begin{align*}
 \abs{ R_{\vt,1} }  = &
 \kappa \Abs{ \intTO{ \Gradedge\vth \cdot \left( \frac{\Grad \tvt}{\tvt} - \Piw \left( \frac{ \Grad \tvt}{\tvt}\right)    \right) }}
 \aleq h \norm{ \Gradedge \vth}_{L^2L^2} \norm{\tvt}_{L^2 W^{2,\infty}},
 \\
\abs{R_{\vt,2}} 
\aleq &\norm{\vth-\tvt}_{L^2L^2}^2 + \norm{ \vthout -\tvt}_{L^2L^2}^2,
\quad
 \abs{R_{\vt,3}} 
\aleq   \frac{1}{\delta}    \norm{\vth-\tvt}_{L^2L^2}^2 + \delta \norm{\Gradedge\vth-\Grad\tvt}_{L^2L^2}^2,
\\
 \abs{R_{\vt,4}} 
 \aleq  & \frac{1}{\delta}    \norm{ \vthout -\tvt}_{L^2L^2}^2 + \delta \norm{\Gradedge\vth-\Grad\tvt}_{L^2L^2}^2,
\quad
\abs{R_{\vt,5}} 
 \aleq h \norm{\Gradedge \vth}_{L^2L^2}^2 \aleq h.
\end{align*}
Moreover, $\norm{ \vthout  - \tvt}_{L^2L^{2}}$ can be controlled as
\begin{align*}
& \norm{ \vthout  - \tvt}_{L^2L^{2}}^2
\leq
\norm{\vth - \tvt}_{L^2L^{2}}^2
+\norm{ \vthout  - \vth^{\rm in}}_{L^2L^{2}}^2
= \norm{\vth - \tvt}_{L^2L^{2}}^2 +h^2\norm{\Gradedge \vth}_{L^2L^{2}}^2 \aleq \norm{\vth - \tvt}_{L^2L^{2}}^2 +h^2.
\end{align*}
Consequently, with \eqref{EN} we have
\begin{equation*}
\abs{ R_{\vt}}
 \aleq h + \int_0^\tau \RE{\vrh,\vuh,\vth}{\tvr,\tvu,\tvt}(t) \dt
+\delta \norm{\Gradedge\vth-\Grad\tvt}_{L^2L^2}^2.
\end{equation*}
Collecting the above estimates, we complete the proof.
\end{document}